\documentclass[12pt]{article}
    \usepackage{srcltx}
   \usepackage{enumerate}
 \usepackage{latexsym}
   \usepackage{indentfirst}
\usepackage{amsmath,amsfonts,amssymb,amsthm}
   \usepackage[russian]{babel}

\newtheorem{theor}{Теорема}
\newtheorem{cor}{Следствие}
\newtheorem{ex}{ Пример}

\newtheorem{lem}{Лемма}

\theoremstyle{remark}
\newtheorem{rem}
{З~а~м~е~ч~а~н~и~е~}

\newcommand{\new}{\newcommand}

\title{
Asymptotic properties of one-step $M$-estimators based on nonidentically distributed
observations.
 }
\author{Yu. Yu. Linke
\thanks{Sobolev Institute of Mathematics,
4 Academician Koptyug ave., Novosibirsk, 630090  Russia;
Novosibirsk State University, 2 Pirogov str., 630090  Russia.
E-mail:
linke@math.nsc.ru}   }
\date{}

\begin{document}
\new{\bb}{\mbox{\mathversion{bold}$\beta$}}
\new{\ee}{\mbox{\mathversion{bold}$\epsilon$}}
\new{\tha}{\mbox{\mathversion{bold}$\theta$}}
\new{\Tha}{\mbox{\mathversion{bold}$\Theta$}}
\new{\et}{\mbox{\mathversion{bold}$\eta$}}
\new{\ophi}{{\overline\varphi}}
\new{\oophi}{{\overline{\overline\varphi}}}
\new{\opsi}{{\overline\psi}}
\new{\pphi}{\mbox{\mathversion{bold}$\varphi$}}
\new{\ppsi}{\mbox{\mathversion{bold}$\psi$}}
\new{\rrho}{\mbox{\mathversion{bold}$\rho$}}
\new{\ddelta}{\mbox{\mathversion{bold}$\delta$}}
\new{\aalpha}{\mbox{\mathversion{bold}$\alpha$}}
\new{\xxi}{\mbox{\mathversion{bold}$\xi$}}

 \maketitle

 {\it Abstract:}
 We study asymptotic behavior of one-step $M$-estimators based on samples from arrays of not necessarily identically distributed random variables and representing
 explicit approximations to the corresponding consistent $M$-estimators. These estimators generalize Fisher's one-step approximations to consistent maximum likelihood estimators. 
As a consequence, we consider some   nonlinear regression problems where the procedure mentioned allow us to construct explicit asymptotically optimal estimators.  We also consider the problem of constructing initial estimators which are  needed for    one-step estimation  procedures.



\section{Введение}

{\bf 1.1.} Пусть  $X_1,\ldots, X_n$ --- независимые не обязательно одинаково распределенные наблюдения произвольной природы, распределения которых зависят от неизвестного параметра $\theta\in \Theta$.
Задача состоит в оценивании этого  параметра по  наблюдениям $X_1,\ldots, X_n$.
  В случае разнораспределенных наблюдений $M$-оценивание принято считать (см., например, [\ref{2007-B}]) одним  из  основных общих методов  получения оценок. В частности, применительно к задачам регрессии этот подход включает в себя метод наименьших квадратов и его модификации, метод максимального правдоподобия,  метод квази-правдоподобия (см., например, [\ref{1997-H}], [\ref{2015-3?}]) и др.

  Определим   $M$-оценки
 как статистики $\widetilde\theta_n$, которые определены с вероятностью, стремящейся к $1$, и являются решениями по $t$ уравнений вида
\begin{equation}\label{p0-0+}
\sum\limits_{i=1}^nM_i(t,X_i)=0
\end{equation}
для некоторого набора функций   $M_i(t,x)$, $i=1,\ldots,n,$ с условием  ${\bf E}M_i(\theta,X_i)=0$ при всех $i$.
Отметим, что далеко не все корни вышеприведенного  $M$-уравнения можно определить на множестве выборок асимптотически полной меры. Так, например, возможна ситуация, когда количество корней при любом $n$ (и даже в пределе при $n\to\infty$, см. [\ref{1985-R}]) будет невырожденной случайной величиной, и лишь один корень можно определить с вероятностью, стремящейся к $1$. Более того,  возможны ситуации, когда несколько корней $M$-уравнения могут быть определены с вероятностью, стремящейся к $1$ (см. пример \ref{p1-prim!}), или таковых корней нет, но с ненулевой вероятностью другие корни могут  существовать.

  Исследованию условий существования и свойств оценок такого типа в различных задачах (в частности, оценок максимального правдоподобия и его различных модификаций, построенных по выборке независимых не обязательно
 одинаково распределенных наблюдений) посвящено больше количество публикаций (см., например, [\ref{2007-B}], [\ref{2004-L1}], [\ref{2004-L2}],   [\ref{1997-H}],  [\ref{1971-H}],  [\ref{1973-P}], [\ref{1975-P}],  [\ref{2004-Y}]).

 C одной стороны, хорошо известно,  что при некоторых  условиях регулярности (см., например, [\ref{2007-B}])
\begin{equation}\label{p0-0-}
\frac{J_n}{\sqrt{I_n}}\big(\widetilde\theta_{n}-\theta\big)\Longrightarrow {\cal N}(0,1),
\end{equation}
где $I_n=\sum\limits_{i=1}^n{\bf E}M_i^2(\theta,X_i)>0, $ $J_n:=\sum\limits_{i=1}^n{\bf E}M_i'(\theta,X_i)\neq 0,$
 через $M_i'(\cdot,\cdot)$ обозначены производные этих функций по первому аргументу, а запись вида
$\eta_n\Rightarrow {\cal N}(0,1)$ означает слабую
сходимость  распределений $\eta_n$  к стандартному  нормальному закону (всюду пределы, если не оговорено противное, берутся при $n\to \infty $). 

 С другой стороны,   поиск  состоятельного  решения $\widetilde\theta_n$ уравнения (\ref{p0-0+}) или поиск его приближения 
  бывает   достаточно сложен, особенно в случае существования большого числа корней уравнения (\ref{p0-0+}) (например, поиск оценки метода наименьших квадратов для задач нелинейной регрессии, связанный с отысканием глобального экстремума функции, см. [\ref{1981-Dem}], [\ref{1989-D}], [\ref{2003-S}]; поиск оценки квази-правдоподобия, см. [\ref{1997-H}]).
 Ситуация существенно упрощается, если известна некоторая  предварительная (состоятельная с требуемой   скоростью сходимости) оценка~$\theta_n^*$ параметра $\theta$.
 В этом случае зачастую  достаточно лишь одного шага  метода Ньютона, чтобы в явном виде построить оценку, обладающую той же асимптотической точностью, как и $M$-оценка $\widetilde\theta_n$, удовлетворяющая (\ref{p0-0-}).

 Основной объект исследования  в первой части  данной   работы --- это статистики $\theta_{n,M}^{**}$,
 имеющие следующую структуру:
\begin{equation}\label{p0-0}
\theta_{n,M}^{**}
=\theta_n^*-{\sum\limits_{i=1}^nM_i(\theta_n^*,X_i)}\big/{\sum\limits_{i=1}^nM_i'(\theta_n^*,X_i)},
\end{equation}
а также некоторые аналоги  этих оценок; здесь и далее $\theta_n^*$ --- некоторая предварительная состоятельная оценка параметра.  Оценка  $\theta_{n,M}^{**}$
  представляет собой одношаговое приближение  корня уравнения (\ref{p0-0+}) методом Ньютона с начальной точкой $t=\theta_n^*$.    В работе будет установлено, что при некоторых предположениях эта оценка  удовлетворяет соотношению
 \begin{equation}\label{p0-0--}
\frac{J_n}{\sqrt{I_n}}\big(\theta_{n,M}^{**}-\theta\big)\Longrightarrow {\cal N}(0,1),
\end{equation}
т.е.  является асимптотически нормальной с той же асимптотической дисперсией, что и $M$-оценка $\widetilde\theta_n$.
  Подчеркнем, что обычно функции $\{M_i(t,x)\}$ в той или иной статистической постановке  выбираются таким образом, чтобы обеспечить желаемые свойства (в том или ином смысле оптимальность) соответствующей $M$-оценки
   (см., например, [\ref{1997-H}], а также примеры в \S3).

Стоит отметить, что в литературе для обозначения таких по  сути  двухшаговых оценок
  используется термин {\it одношаговые} оценки в том смысле, что
при наличии некоторой {\it предварительной оценки} (существование которой зачастую лишь постулируется) именно  {\it за один шаг} метода Ньютона  происходит качественное  улучшение этой начальной статистики.
Следуя  общепринятой в последние десятилетия  терминологии (начиная с работ  П. Хьюбера,
П. Бикела, Р. Серфлинга, Н. Веравербеке,  П.-К. Сена, К. Мюллер,   Х. Юречковой,  Д. Рупперта, Р. Кэрролла, Дж. Фана  и др.)
оценки такого типа будем называть {\it одношаговыми} $M$-оценками.

{\bf 1.2.} В этом разделе приведем краткий обзор известных автору результатов об одношаговых $M$-оценках.

Впервые  идея одношагового  оценивания была предложена, по-видимому,   Р.~Фишером [\ref{1925-F}].  В качестве приближений для состоятельной оценки
максимального правдоподобия $\widetilde \theta_n$, построенной по
 выборке с  плотностью   $f(\theta,x)  $
    относительно некоторой $\sigma$-конечной  меры,  Р. Фишер предложил использовать
одну из следующих двух статистик:
\begin{equation}\label{2}
\theta_{n,L}^{**}
=\theta_n^*-\frac{L'_n(\theta_n^*)}{L''_n(\theta_n^*)}\
\qquad\mbox{или}\qquad \theta_{n,I}^{**}
=\theta_n^*+\frac{L'_n(\theta_n^*)}{nI(\theta_n^*)},
\end{equation}
где
$I(t)$ --- информация Фишера, соответствующая плотности
$f(t,x),$  $L_n(t)$ --- логарифмическая функция правдоподобия,
построенная по   первым $n$ элементам выборки, т.е.
$L_n(t)=\sum\limits_{i=1 }^n\ln f(t,X_i)$, а $\theta_n^*$ --- нередко оценка метода моментов.
Статистики (\ref{2}) для краткости в дальнейшем будем  называть  {\it оценками Фишера}.

Выделим прежде всего исследования, касающиеся одношаговых оценок в случае одинаково распределенной выборки.
 Ряд монографий и статей  (см., например, [\ref{2007-B}],  [\ref{1975-Z}], [\ref{1985-J}], [\ref{1990-J}], [\ref{1991-L}], [\ref{1984-H}], [\ref{1956-L}], [\ref{1980-S}], [\ref{2007-V}], [\ref{2014-1}], [\ref{2015-1?}], [\ref{2222}]  и ссылки там же) содержат исследования оценок Фишера,  либо их обобщений на случай приближенного поиска $M$-оценок  в случае одинаково распределенной выборки  (в этом случае функции $M_i(\cdot,\cdot)$ в (\ref{p0-0}) не зависят от $i$).
Достаточные условия асимптотической  нормальности
одношаговых оценок в этом частном случае
 найдены в [\ref{1980-S}], 
 [\ref{1975-Z}], [\ref{1991-L}], [\ref{1956-L}], [\ref{2014-1}], [\ref{2222}].
 Особо выделим  работу Л. Ле Кама   [\ref{1956-L}], в которой, по-видимому, впервые отмечено, что  полученные в указанной работе  условия асимптотической нормальности одной из  оценок Фишера не только не влекут состоятельность  оценки максимального правдоподобия, но даже не гарантируют существование этой оценки (см. также     [\ref{1991-L}], [\ref{2014-1}], [\ref{2222}] и пример \ref{p1-pr-L}).
Кроме того, помимо оценивания  статистиками вида  $\theta_{n,M}^{**}$
параметра~$\theta$,
зачастую изучается  близость этих  же  оценок    с состоятельной $M$-оценкой $\widetilde\theta_n$.  Так, в [\ref{2007-B}]   доказано (см. также ссылки в
 [\ref{1985-J}]), что
\begin{equation}\label{p0-4}
\sqrt{n}\big(\theta_{n,M}^{**}-\widetilde\theta_n\big)\stackrel{p}
{\to}0
\end{equation}
в случае, когда $\widetilde\theta_n$ --- оценка максимального правдоподобия.
   В  [\ref{1985-J}] и
[\ref{2007-V}] установлено, что  скорость сходимости в (\ref{p0-4}) при некоторых дополнительных предположениях
может быть улучшена и, кроме того, там  доказано, что
${n}\big(\theta_{n,M}^{**}-\widetilde\theta_n\big)=O_p(1)$.
В [\ref{1990-J}]
  найдено предельное
распределение  для
$n\big(\theta_{n,M}^{**}-\widetilde\theta_n\big)$ (не являющееся
нормальным).

Подчеркнем,   что практически  во всех известных автору
работах  в качестве предварительной оценки
 $\theta_n^*$  рассматриваются {\it $\sqrt{n}$-ограниченные} оценки $\theta_n^*$, т.е.
 оценки, удовлетворяющие соотношению
$\sqrt{n}(\theta_n^*-\theta)=O_p(1)$.
Исключение составляют  монография  [\ref{1975-Z}] и работы [\ref{2014-1}] и [\ref{2222}]. В [\ref{1975-Z}]
(теорема 5.5.4)  при  доказательстве асимптотической нормальности
 оценки Фишера $\theta_{n,L}^{**}$
   предполагается  выполненным
  условие   ${n}^{1/4}$-состоятельности  оценки $\theta_n^*$, т.е. сходимость
  ${n}^{1/4}(\theta_n^*-\theta)\stackrel{p}{\to}0$.
 В   [\ref{2014-1}] и [\ref{2222}]  асимптотическая нормальность оценок Фишера
 и одношаговых $M$-оценок   доказана при  более широком спектре ограничений на точность
 предварительной оценки --
начиная от $\sqrt{n}$-ограниченных оценок $\theta_n^*$ и заканчивая
 ${n}^{1/4}$-состоятельными оценками. 
 Там же установлено, что условие ${n}^{1/4}$-состоятельности предварительной оценки
 по сути  необходимо для доказательства асимптотической нормальности оценок Фишера и одношаговых
 $M$-оценок.

 Задача приближенного поиска состоятельных $M$-оценок возникает и в случае, когда скорость сближения предварительной оценки и  параметра оказывается  медленнее, чем $n^{-1/4}$  (например, в задаче оценивания параметра однопараметрических семейств распределений с так называемыми тяжелыми хвостами; см. [\ref{2015-1?}]).
Если предварительная оценка не достаточно точна для непосредственного использования  одношаговой процедуры  улучшения, то точность предварительной оценки может  быть улучшена многократным   применением метода Ньютона-Фишера. Подобная процедура хорошо известна (см., например, [\ref{2007-C}], [\ref{1988-R}]).
Альтернативный подход для однопараметрических семейств распределений предложен автором в [\ref{2015-1?}], где    приведен алгоритм построения новых одношаговых оценок, специально ориентированных  на медленно сближающиеся с параметром предварительные оценки. Эти новые одношаговые оценки отличаются от оценок Фишера (\ref{2}) наличием некоторых дополнительных поправочных слагаемых (зависящих от скорости сближения предварительной оценки и  параметра) и  позволяют за одну итерацию
$n^{\beta}$-состоятельную  оценку при $\beta< 1/4$  улучшать до асимптотически эффективной.

Все, сказанное выше, относилось к случаю однородной выборки. Возможность использования методологии одношагового оценивания в статистических задачах с разнораспределенными выборочными наблюдениями хорошо  известна (см., например, монографии [\ref{2007-B}], [\ref{1984-H}], [\ref{2006-J}], [\ref{2012-J1}], [\ref{1997-H}], [\ref{1989-M}]).
 По-видимому, впервые одношаговые оценки в случае разнораспределенных наблюдений были предложены в
[\ref{1975-B}] для специальных оценок в линейных моделях. Эти одношаговые оценки являются приближениями для  оценок робастного варианта метода наименьших квадратов в задачах линейной регрессии, предложенных в  [\ref{1973-H}].
  В [\ref{1975-B}] доказана асимптотическая нормальность одношаговых оценок в этой специальной задаче.  В [\ref{1987-J}], [\ref{2002-W}] изучаются другие аспекты   одношаговых оценок для задач робастного оценивания в линейных моделях.
 В частности, в [\ref{1987-J}] доказано, что если $n^\tau|\theta_n^*-\theta|=O_p(1)$ при  $\tau\in(1/4,1/2]$,  то $\theta_n^{**}-\widetilde\theta_n=O_p(n^{-2\tau})$, где $\theta_n^{**}$ --- одношаговая статистика из  [\ref{1975-B}], а $\widetilde\theta_n$ --- $M$-оценка из [\ref{1973-H}].

В случае разнораспределенной  (и даже не обязательно состоящей из независимых компонент) выборки методология  одношагового оценивания успешно   используется  в различных специальных статистических задачах
  (см., например,  [\ref{1975-B}], [\ref{1984-H}],
 [\ref{1987-J}] [\ref{2008-Z}],
 [\ref{2000-C}], [\ref{1999-F1}],  [\ref{1999-F}],  [\ref{2006-J}], [\ref{2012-J1}], [\ref{2012-J2}], [\ref{1994-M1}], [\ref{1994-M2}], [\ref{1986-S}], [\ref{1992-S}] и ссылки там же).
Стоит отметить,  что и в случае разнораспределенной выборки  в качестве предварительных оценок параметра используются, как правило, $\sqrt{n}$-ограниченные оценки (точнее, при исследовании одношаговых оценок существование $\sqrt{n}$-ограниченной  предварительной оценки чаще всего просто постулируется).

Исследования асимптотического поведения одношаговых $M$-оценок
 вида (\ref{p0-0})  без той или иной спецификации функций $\{M_i(\cdot,X_i)\}$
  в  случае разнораспределенной выборки  автору не известны.

{\bf 1.3.} О структуре работы.  В \S 2  в широком спектре ограничений на точность предварительной оценки $\theta_n^*$ найдены (теорема \ref{p1-t1}) достаточно общие  условия
асимптотической нормальности (\ref{p0-0--}) одношаговых $M$-оценок,  определяемых соотношением  (\ref{p0-0}).
 При этом
получено
универсальное условие, связывающее
  гладкость функций $\{M_i(t,x)\}$, определяющих
  одношаговые $M$-оценки,   и  точность
  предварительной оценки $\theta_n^*$,   которые
нужны для  асимптотической нормальности этих оценок.
Кроме того, обсуждается вопрос  о  минимальном достаточном ограничении на точность предварительной оценки  $\theta_n^*$, возникающем при доказательстве асимптотической нормальности одношаговых $M$-оценок (теорема \ref{p1-t1w+}).
Построены оценки для асимптотической дисперсии одношаговых  $M$-оценок    и доказана сходимость к стандартному нормальному закону погрешностей оценивания при замене асимптотических дисперсий их оценками (теорема \ref{p1-t2}),  что может быть полезным при построении доверительных интервалов и для проверки статистических  гипотез.
В разделе 2.2 исследован  некоторый аналог  статистики $\theta_{n,M}^{**}$, обобщающий     оценку Фишера $\theta_{n,I}^{**}$ из (\ref{2}) на случай разнораспределенной выборки (теорема \ref{p1-t1!-!!}).
  В разделе 3.3 показано  как можно ослабить условия на гладкость функций $\{M_i(t,X_i)\}$, определяющих одношаговые $M$-оценки, в случае специального вида   этих  функций.

Подчеркнем, что полученные в работе общие условия асимптотической нормальности одношаговых $M$-оценок не гарантируют существование  состоятельной  $M$-оценки.
Таким образом, одношаговые $M$-оценки могут обладать некоторыми нужными свойствами {\it вне зависимости от поведения $M$-оценок}! Т.е. и в случае разнораспределенной выборки имеет место  эффект, отмеченный ранее   Л.Ле Камом [\ref{1956-L}] в частном случае оценок максимального правдоподобия для однородных выборок. Но чем обусловлен такой неожиданный эффект?
И  что же все-таки в этом случае приближают одношаговые $M$-оценки, являясь (по построению) приближением для решения уравнения (\ref{p0-0+})? В разделе 2.3  (теорема \ref{p1-t3})
приведено поясняющее указанную ситуацию утверждение  (см. также пример \ref{p1-prim!}).
 Дело в том, что при некоторых
 весьма слабых ограничениях (не гарантирующих  существование состоятельной $M$-оценки) в окрестности  $\theta$
с вероятностью, стремящейся к~$1$, существует  единственный  корень $\widetilde\theta_n(\theta)$  уравнения~(\ref{p0-0+}) (ближайший к $\theta$ корень этого уравнения), обладающий требуемой точностью.
  Указанную
  случайную величину $\widetilde\theta_n(\theta)$, зависящую от $\theta$ и не являющуюся статистикой,   и приближают одношаговые $M$-оценки. Указанный факт   влечет    асимптотическую нормальность одношаговых   $M$-оценок вне зависимости от поведения $M$-оценки: лишь бы существовал корень $\widetilde\theta_n(\theta)$, обладающий нужными свойствами.

Таким образом,  одношаговые $M$-оценки   в случае
разнораспределенной выборки могут быть асимптотически нормальными и
без предположений о том или ином поведении собственно самих
$M$-оценок. При этом,  если существует состоятельная $M$-оценка
$\widetilde\theta_n$, то с вероятностью,
стремящейся к $1$, оценка $\widetilde\theta_n$ и величина
$\widetilde\theta_n(\theta)$ совпадают и одношаговые оценки
являются  приближениями для состоятельных $M$-оценок. Отметим еще,
что ранее понятие ближайшего к~$\theta$ корня  уравнения
правдоподобия упоминалось   в  [\ref{1991-L}] и [\ref{2015-1?}] в
случае одинаково распределенной выборки, а поведение ближайшего к
$\theta$ корня $\widetilde\theta_n(\theta)$ в случае одинаково
распределенной выборки исследовалось в [\ref{2222}]. Тот факт, что и
в случае разнораспределенных наблюдений одношаговые $M$-оценки могут
приближать локальный экстремум соответствующей функции, а не глобальный, отмечается,
например,  [\ref{1975-B}], [\ref{2014-F}] (имеются в виду ситуации, когда $M$-оценки порождаются некоторой экстремальной задачей).

  Наш интерес к одношаговым $M$-оценкам в случае разнораспределенных наблюдений в первую очередь связан с возможностью эффективного применения подобных оценок  в различных регрессионных задачах (см. также [\ref{2015-3?}]).
  К важному базовому   предположению  методологии одношагового  оценивания  
  можно отнести наличие некоторой  достаточно хорошей предварительной оценки.
В разделе 3.1   предлагаемая процедура одношагового  улучшения иллюстрируется  несколькими  примерами   моделей нелинейной регрессии, для которых известны  явные, но не оптимальные оценки.  В разделе 3.2  предлагается  некоторый  способ построения   предварительных оценок в задачах нелинейной регрессии.

Обобщения основных  результатов \S2  на случай многомерного основного параметра содержатся в  \S4.
 Доказательства всех утверждений отнесены в~\S5. %

\section{Асимптотические свойства одношаговых $M$-оценок}

{\bf 2.1.} Нам потребуются следующие условия.

$ ({\bf A}_1).$
Наблюдается выборка объема $n$,  состоящая из независимых элементов   $X_{1},$ $X_{2},\ldots, X_{n}$  со значениями в
 произвольном измеримом пространстве
${\cal X}$ и  распределениями~${\cal L}_{1,\theta},$ ${\cal L}_{2,\theta}\ldots,{\cal L}_{n,\theta}$, зависящими  от    интересующего нас  основного
неизвестного параметра $\theta\in\Theta\subset
\mathbb{R}$, где $\Theta$ --- открытое множество.  Кроме того,  эти распределения зависят, вообще говоря, от $n$, и, возможно, зависят  также  от некоторых мешающих  параметров
 произвольной природы.

$ ({\bf A}_2).$  На множестве  $\Theta$ при любом  $i$ заданы зависящие, вообще говоря, от $n$  функции  $M_i(t,X_i)$ и $M'_i(t,X_i)$
такие, что
для каждого интервала
$(t_1,t_2)$, целиком лежащего в $\Theta$, с
   вероятностью $1$ имеет место равенство
$$
M_i(t_2,X_i)-
M_i(t_1,X_i)=\int\limits_{t_1}^{t_2}
 M'_i(t,X_i)dt\qquad \forall (t_1,t_2)\subset  \Theta
 $$
 и справедливы  соотношения
\begin{equation}\label{p1-25}
{\bf E}M_i(\theta,X_i)=0, \qquad {\bf E}M^2_i(\theta,X_i)<\infty, \qquad
{\bf E}|M'_i(\theta,X_i)|<\infty.
\end{equation}
Здесь и всюду в дальнейшем символ $\theta$ обозначает истинное значение основного параметра распределений  и
   математическое ожидание в (\ref{p1-25}) и далее  берется по
    распределению ${\cal L}_{i,\theta}$. Отметим еще, что  для почти всех $t\in\Theta$ функции $M'_i(t,X_i)$ с вероятностью $1$ совпадают с частной производной функции $M_i(t,X_i)$ по первому аргументу.

$ ({\bf A}_3).$ Начиная с некоторого $n$,
$$
I_n:=
\sum\limits_{i=1}^n
{\bf E}M^2_i(\theta,X_i)>0, \qquad J_n:=
\sum\limits_{i=1}^n{\bf E}M_i'(\theta,X_i)\neq 0
$$
и, кроме того,
$$\quad |J_n|/\sqrt{I_n}\to \infty,$$
\begin{equation}\label{p1-5}
\sum\limits_{i=1}^n M'_i(\theta,X_i)/J_n\stackrel{p}{\to}1, \qquad
\sum\limits_{i=1}^n M_i(\theta,X_i)/\sqrt{I_n}\Longrightarrow {\cal N}(0,1).
\end{equation}

Положим $\quad \overline\tau_n(\delta)=|J_n|^{-1}\displaystyle\sum\limits_{i=1}^n{\bf E}\tau_{i}(\delta,X_i),$ где
 \begin{equation}   \label{p1-tau}
\tau_{i}(\delta,X_i):=
\begin{cases}
\sup\limits_{t
:\; |t-\theta|\leq\delta}\big|
M'_i(t,X_i)-M'_i(\theta,X_i) \big|,\quad \mbox{если }\; [\theta-\delta,\theta+\delta]\subset\Theta,
\\
\; \infty,\quad \mbox{иначе}.
\end{cases}
\end{equation}

$ ({\bf A}_4).$   Имеет место соотношение $\lim\limits_{\delta\to 0}\limsup\limits_{n\to\infty} \overline\tau_n(\delta)=0.
$

$ ({\bf A}_5).$ Имеется  оценка  $\theta_n^*=\theta_n^*(X_1,\ldots,X_n)$ такая, что
 \begin{gather*}
  \frac{|J_n|}{\sqrt{I_n}}|\theta_n^*-\theta|\overline\tau_n(|\theta_n^*-\theta|)\stackrel{p}{\to} 0.
  \end{gather*}

Группу условий, состоящую из
   предположений $(A_1),\ldots, $   $(A_5)$,  будем далее  обозначать  символом $(A) $.

При выполнении условий  $(A) $ в качестве оценки параметра $\theta$ рассмотрим   статистику $\theta_{n,M}^{**}$, определяемую   соотношением (\ref{p0-0}), т.е.
 \begin{gather} \label{p1-7}
 \theta_{n,M}^{**}=\theta_n^*- {\sum\limits_{i=1}^n M_i(\theta_n^*,X_i)}\big/{
 \sum\limits_{i=1}^nM'_i(\theta_n^*,X_i)}
 \end{gather}
 во всех случаях, когда определены все величины в правой  части  равенства в~(\ref{p1-7}).
Тем самым, допускается, что знаменатель отношения в правой  части  равенства~(\ref{p1-7})
 может обращаться в ноль, или что функции в правой части этого равенства не определены.
Но, как будет установлено при доказательстве  теоремы~\ref{p1-t1},   при выполнении
  условия $(A) $ справедливы соотношения  ${\bf P}(\theta_n^*\in\Theta)\to 1$ и ${\bf P}\Big(\sum\limits_{i=1}^nM'_i(\theta_n^*,X_i)\neq 0\Big)\to 1,$ а потому
статистика
 $\theta_{n,M}^{**}$  определена с вероятностью, стремящейся к $1$.

Имеют место следующие два утверждения.
\begin{theor}                                     \label{p1-t1}
Пусть  выполнены  условия   $(A)$.
 Тогда
оценка $\theta_{n,M}^{**}$  определена с вероятностью, стремящейся к $1$ и
\begin{gather}                                                                                           \label{p1-1w}
\Delta_{n,M}:=\frac{J_n}{\sqrt{I_n}}\big(\theta_{n,M}^{**}-\theta\big)+\sum\limits_{i=1}^nM_i(\theta,X_i)/\sqrt{I_n}\stackrel{p}{\to}0.
 \end{gather}
В частности,
\begin{gather}                                                                                           \label{p1-24}
\frac{J_n}{\sqrt{I_n}}\big(\theta_{n,M}^{**}-\theta\big)\Longrightarrow {\cal N}(0,1).
 \end{gather}
\end{theor}

 \begin{theor}                                     \label{p1-t2}
Если   выполнены  условия   $(A)$ и
\begin{equation} \label{p1-8}
\sum\limits_{i=1}^n [M_i'(\theta,X_i)]^2/J_n^2\stackrel{p}{\to}0,\quad\mbox{}\quad
\sum\limits_{i=1}^n M_i^2(\theta,X_i)/I_n\stackrel{p}{\to}1,
\end{equation}
то
\begin{equation} \label{p1-9}
d_{n,M}^*\big(\theta_{n,M}^{**}-\theta\big)\Longrightarrow {\cal N}(0,1)\quad\mbox{при}\quad d_{n,M}^*=\sum\limits_{i=1}^n M_i'(\theta_n^*,X_i)\big/
\Big(\sum\limits_{i=1}^n M_i^2(\theta_{n,M}^{**},X_i)\Big)^{1/2}
\end{equation}
и статистика $d_{n,M}^*$  определена с вероятностью, стремящейся к $1$.
\end{theor}
Теорема \ref{p1-t2}   может быть полезной при построении доверительных интервалов и проверке гипотез, поскольку нормирующий множитель  $d_{n,M}^*$ у разности   $\theta^{**}_{n,M}-\theta$ является статистикой.

\begin{rem} \label{p1-z1}{\rm
Соотношения (\ref{p1-5}) в условии $(A_3)$ --- это соответственно  варианты закона больших чисел и центральной предельной теоремы в схеме серий. Здесь для их выполнения
достаточно потребовать следующие условия (см., например, [\ref{2007-B}]):
$$
\sum\limits_{i=1}^n{\bf E}\min \left\{\big|M_i'(\theta,X_i)-{\bf E}M_i'(\theta,X_i)\big|/|J_n|,\; \big|M_i'(\theta,X_i)-{\bf E}M_i'(\theta,X_i)\big|^2/J_n^2\right\}\to 0,
$$
$$
\sum\limits_{i=1}^n{\bf E}\min \left\{\big|M_i^2(\theta,X_i)\big|/I_n,\; \big|M_i(\theta,X_i)\big|^s/I_n^{s/2}\right\}\to 0
$$
при некотором $s>2.$ Аналогично можно выписать достаточное условие и для выполнения закона больших чисел, возникающего в теореме \ref{p1-t2} (см. второе условие в (\ref{p1-8})).
Отметим еще, что в случае разнораспределенных наблюдений  первого соотношения в (\ref{p1-5}) не достаточно, чтобы было выполнено первое  условие в (\ref{p1-8}). Соответствующий пример нетрудно построить.

 Отметим еще, что условие независимости из $(A_1)$ в теореме \ref{p1-t1} может быть ослаблено. При отказе от условия независимости числовую характеристику $I_n$ нужно определить равенством  $I_n\equiv {\bf D}\Big(\sum\limits_{i=1}^n M_i(\theta,X_i)\Big)$. Достаточные  условия для выполнения  соотношений (\ref{p1-5}) для тех или иных видов зависимостей также хорошо известны.~\hfill$\square$
 }
\end{rem}

\begin{rem} \label{p1-zz1}
Обычно предполагается, что $\Theta$ -- открытое множество, но
 такого ограничения нет  в работе Л. Ле Кама [\ref{1956-L}], где доказано
утверждение об асимптотической нормальности оценок Фишера для однородной выборки.
Чтобы   теорема \ref{p1-t1} об асимптотической нормальности одношаговых $M$-оценок  обобщала соответствующий результат из [\ref{1956-L}] в случае одномерного параметра,  достаточно в $(A_1)$  считать, что $\Theta$ --- произвольное  множество, а  в условиях $(A_2)-(A_5)$ множество  $\Theta$ заменить на  некоторое  открытое множество $\Theta_o\supset\Theta$. В этом случае все приведенные выше утверждения  сохранятся.

Следствия из  теоремы  \ref{p1-t1} и приводимой далее теоремы \ref{p1-t1!-!!}  на случай однородной выборки  обобщают также  соответствующие результаты из работ [\ref{1980-S}], [\ref{1975-Z}], [\ref{1991-L}], [\ref{2014-1}] 
об  асимптотической нормальности оценок Фишера или  одношаговых $M$-оценок  в случае одинаково распределенной выборки.
\hfill$\square$
\end{rem}

\begin{rem}\label{p1-z15}
Центральное условие группы предположений $(A)$ --- это условие $(A_5)$.
 Это условие является  некоторым универсальным ограничением, связывающем гладкость функций $\{M_i(\cdot,X_i)\}$, определяющих одношаговую $M$-оценку $\theta_{n,M}^{**}$, со скоростью сближения предварительной оценки $\theta_n^*$  и параметра~$\theta$,
    которые
нужны для  асимптотической нормальности одношаговых $M$-оценок. При этом точность предварительной оценки $\theta_n^*$ и гладкость функций
$\{M_i(\cdot,X_i)\}$ обратно пропорциональны друг другу:  чем меньше точность, тем больше должна быть гладкость.
Действительно,  для выполнения указанного ограничения   достаточно, чтобы
 величина
 $\overline\tau_n(\delta)$  и оценка $\theta_n^*$ одновременно (для некоторого  $\alpha$ из указанной области) удовлетворяли условиям одной из следующих двух групп:
\begin{equation} \label{p1-415}
\limsup\limits_{n\to\infty} \overline\tau_n(\delta)=o(\delta^\alpha)\quad\mbox{при}\quad\delta\to 0\quad\mbox{и}\quad \Big(\frac{|J_n|}{\sqrt{I_n}}\Big)^{1/(1+\alpha)}|\theta_n^*-\theta|=O_p(1),\quad 0\leq \alpha<1;
\end{equation}
\begin{equation}\label{p1-416}
\limsup\limits_{n\to\infty}\overline\tau_n(\delta)=O(\delta^\alpha)\quad\mbox{при}\quad\delta\to 0\quad\mbox{и}\quad \Big(\frac{|J_n|}{\sqrt{I_n}}\Big)^{1/(1+\alpha)}|\theta_n^*-\theta|\stackrel{p}{\to}0,\quad 0< \alpha\leq 1.
\end{equation}
В частности, в двух крайних случаях $\alpha=0$ и $\alpha=1$ получаем, что  для справедливости  $(A_5)$ достаточно, чтобы
$$\limsup\limits_{n\to\infty} \overline\tau_n(\delta)=o(1)\quad\mbox{при}\quad\delta\to 0\qquad \mbox{и}\qquad\frac{|J_n|}{\sqrt{I_n}}|\theta_n^*-\theta|=O_p(1), $$
либо
$$
 \limsup\limits_{n\to\infty}\overline\tau_n(\delta)=O(\delta)\quad\mbox{при}\quad\delta\to 0\qquad \mbox{и}\qquad \frac{|J_n|^{1/2}}{I_n^{1/4}}|\theta_n^*-\theta|\stackrel{p}{\to}0.
$$
Как будет установлено  далее (см. теорему \ref{p1-t1w+}),  при широких ограничениях    условие
 \begin{equation}                                                                   \label{p1-5w}
\frac{|J_n|^{1/2}}{I_n^{1/4}}|\theta_n^*-\theta|\stackrel{p}{\to}0
 \end{equation}
по сути необходимо для доказательства асимптотической нормальности одношаговых $M$-оценок $\theta_{n,M}^{**}.$
Таким образом, асимптотическая нормальность одношаговых оценок установлена в достаточно широком  спектре ограничений на точность предварительной оценки $\theta_n^*$.
Отметим еще, что если $|J_n|/\sqrt{I_n}$ имеет порядок $\sqrt{n},$ то условия на точность оценки  $\theta_n^*$ в указанных  двух крайних случаях --- это, соответственно, либо предположение о $\sqrt{n}$-ограниченности, либо об $n^{1/4}$-состоятельности этой оценки.
 \hfill$\square$
\end{rem}
\begin{rem}
Помимо различных ``регулярных'' вариантов поведения двух рассматриваемых в замечании \ref{p1-z15} характеристик, нельзя исключать и, в известной степени,  вырожденный случай,  когда c вероятностью 1 функции  $M_i'(t,X_i)$ постоянны в некоторой
  окрестности $
  \{t
 : |t-\theta|\leq \delta(\theta)\}$     точки~$\theta$, т.~е. $\tau_i(\delta,X_i)=\overline\tau_n(\delta)\equiv 0$ при $\delta\leq \delta(\theta)$. В этом случае величина~$\theta_n^*$ может быть  произвольной и даже не сходиться к $\theta$, лишь бы $\theta_n^*$
 попадала в указанную окрестность $\theta$ с вероятностью, стремящейся к $1$.  Отметим, что подобная ситуация реализуется для оценки максимального правдоподобия параметра сноса   $\theta$  для однородной выборки из нормального распределения.\hfill$\square$
\end{rem}

Обсудим подробнее  условия  на точность предварительной оценки~$\theta_n^*$, используемые в теореме \ref{p1-t1}. Оказывается, при некоторых
  достаточно широких ограничениях сходимость (\ref{p1-5w})
  является необходимым условием для справедливости (\ref{p1-1w}).
  Таким образом, условие  (\ref{p1-5w})
можно считать в некотором смысле минимальным достаточным ограничением для доказательства
  (\ref{p1-24}), поскольку    метод доказательства,  используемый для вывода асимптотической нормальности  одношаговых оценок как в данной работе, так и других исследованиях  (см.  все указанные в п.1.2 ссылки)
требует от предварительной оценки $\theta_n^*$ скорость сходимости не хуже, чем в условии (\ref{p1-5w}).
  Чтобы сформулировать соответствующее утверждение, нам потребуется дополнительное предположение.

$(A_6)$ При всех $i$
функции  $M_i(t,X_i)$ c вероятностью $1$ дважды  дифференцируемы  по $t$
и для любого $\varepsilon>0$
$$P(\delta):=\limsup\limits_{n\to\infty}{\bf P}\Big(\sum\limits_{i=1}^n\sup\limits_{t\in \Theta:|t-\theta|\leq \delta}|M''_i(t,X_i)-M''_i(\theta,X_i)|/|J_n|>\varepsilon\Big)\to 0 \quad \mbox{при}\;\delta\to 0. $$
 Кроме того,
$\overline{M_n''}(\theta):=\sum\limits_{i=1}^n\big(M_i''(\theta,X_i)-{\bf E}M_i''(\theta,X_i)\big)/J_n\stackrel{p}{\to}0.$

Справедлива
\begin{theor}                                     \label{p1-t1w+}
Пусть выполнены     условия  $(A)$,   $(A_6)$ и    предварительная оценка   $\theta_n^*$ состоятельна.
 Тогда   справедливо соотношение
\begin{equation}                                                                                                 \label{p1-3w}
\Delta_{n,M}= \frac{J_n}{\sqrt{I_n}}(\theta_n^*-\theta)^2
\Big(\sum\limits_{i=1}^n {\bf E}M_i''(\theta,X_i)/(2J_n)+o_p(1)\Big)+o_p(1).
 \end{equation}

 Если дополнительно
$\liminf\limits_{n\to\infty}\big|\sum\limits_{i=1}^n {\bf E}{M}_i''(\theta,X_i)\big|/|J_n|\geq \delta$ при некотором $\delta>~0$,
 то  условие {\rm (\ref{p1-5w})}
 является
необходимым и
достаточным для того, чтобы сходимость {\rm(\ref{p1-1w})} имела место.
 \end{theor}

{\bf 2.2.} Рассмотрим некоторую модификацию оценки $\theta_{n,M}^{**}, $ определенной в  (\ref{p1-7}).
Поскольку при широких ограничениях $J_n\equiv \sum\limits_{i=1}^n {\bf E}M_i'(\theta,X_i)$  в некотором смысле приближает величину $\sum\limits_{i=1}^n M_i'(\theta_n^*,X_i),$ то  в качестве одношаговой  оценки параметра $\theta$ можно рассматривать
 оценку $ \theta_{n,J}^{**},$ определяемую соотношением
 \begin{gather} \label{p1-777}
 \theta_{n,J}^{**}=\theta_n^*- {\sum\limits_{i=1}^n M_i(\theta_n^*,X_i)}\big/{
 J_n^*},
 \end{gather}
где $J_n^*$ --- некоторая статистика, приближающая $J_n.$
 В частном случае  однородной выборки подобная  модификация одного шага  метода Ньютона была  предложена Р.Фишером в задаче построения одношаговых приближений  для состоятельных оценок максимального правдоподобия (см. (\ref{2})).
 В дальнейшем подобный прием получил название ``method of scoring'' (см., например, [\ref{1992-O}]).

 В случае разнораспределенных наблюдений одношаговая оценка  вида (\ref{p1-777})
 впервые была предложена, по-видимому, в [\ref{1975-B}] для одной специальной статистической задачи.
 В некоторых задачах выбор статистики  $J_n^*$ может быть достаточно прост.  Например, в случае построения одношаговых приближений для оценки максимального правдоподобия величина  $-J_n=-J_n(\theta)$ есть суммарная информация Фишера, построенная по элементам разнораспределенной выборки. Таким образом,  можно положить  $J_n^*=J_n(\theta_n^*)$.

    В общем случае распределения наблюдений могут зависеть от некоторых мешающих параметров и
  выбор той или иной  оценки $J_n^*$ может быть осуществлен только в рамках конкретной статистической модели (см. примеры в \S3).
    Поэтому  аналог утверждений  теорем \ref{p1-t1} и \ref{p1-t2} приведем в предположении, что существует статистика  $J_n^*$ такая, что $J_n^*/J_n\stackrel{p}{\to}1$.

\begin{theor}                                     \label{p1-t1!-!!}
Пусть  выполнены условия   $(A)$ и
\begin{gather}                                                                                           \label{p1-555-}
\frac{|J_n|}{\sqrt{I_n}}|\theta_n^*-\theta|\Big(|J_n^*/J_n-1|+\big|\overline {M_n'}(\theta)\big|\Big)\stackrel{p}{\to}0, \qquad J_n^*/J_n\stackrel{p}{\to}1,
 \end{gather}
 где  $\overline {M_n'}(\theta):=\sum\limits_{i=1}^n M_i'(\theta,X_i)/J_n-1.$
 Тогда одношаговая $M$-оценка $\theta_{n,J}^{**}$  определена с вероятностью, стремящейся к~$1$~и
\begin{gather}                                                                                           \label{p1-24g-}
\frac{J_n}{\sqrt{I_n}}\big(\theta_{n,J}^{**}-\theta\big)\Longrightarrow {\cal N}(0,1).
 \end{gather}
 Если дополнительно справедливо  {\rm(\ref{p1-8})}, то
\begin{equation} \label{p1-9g-}
d_{n,J}^*\big(\theta_{n,J}^{**}-\theta\big)\Longrightarrow {\cal N}(0,1)\quad\mbox{при}\quad d_{n,J}^*=J_n^*\big/
\Big(\sum\limits_{i=1}^n M_i^2(\theta_{n,J}^{**},X_i)\Big)^{1/2}
\end{equation}
и статистика $d_{n,J}^*$  определена с вероятностью, стремящейся к $1$.
\end{theor}

Таким образом, асимптотические свойства  одношаговой $M$-оценки $\theta_{n,J}^{**}$ аналогичны поведению  одношаговой $M$-оценки $\theta_{n,M}^{**}$.
Более того, за счет возможного многообразия различных приближений   $J_n^*$ для $J_n$ оценка
  $\theta_{n,J}^{**}$ может определяться не однозначно (см. пример \ref{p1-e2}). Тем самым, в некоторых задачах можно предложить множество явных одношаговых оценок,  эквивалентных в смысле асимптотической точности. Понятно, что среди этих оценок в тех или иных ситуациях можно выбирать  в том или ином смысле более предпочтительную оценку для параметра $\theta$.

\begin{rem}\label{p1-z15+}
Пусть $J_n^*=J_n(\theta_n^*)$ (и, следовательно, $J_n(\theta)=J_n$). Приведем в этом случае простые достаточные условия для выполнения условия   (\ref{p1-555-}).
 В \S5 будет установлено, что для выполнения $(A_5)$ и  условия   (\ref{p1-555-}) (здесь считаем, что имеют место условия $(A_1)$-$(A_4)$) достаточно, чтобы величины  $\overline \tau_n(\delta),$ $J_n(t)$    и оценка $\theta_n^*$ одновременно удовлетворяли одному из следующих двух блоков  условий:

1) выполнено условие (\ref{p1-415}) и
\begin{equation} \label{p1-415!}
 \limsup\limits_{n\to\infty}\Big|\frac{J_n(t)}{J_n}-1\Big|=o(|t-\theta|^\alpha)\quad\mbox{при}\quad t\to \theta\quad\mbox{и}\quad \Big(\frac{|J_n|}{\sqrt{I_n}}\Big)^{\alpha/(1+\alpha)}\big|\overline {M_n'}(\theta)\big|\stackrel{p}{\to }0;
\end{equation}

2) выполнено условие (\ref{p1-416}) и
\begin{equation}\label{p1-416!}
 \limsup\limits_{n\to\infty}\Big|\frac{J_n(t)}{J_n}-1\Big|=O(|t-\theta|^\alpha)\;\mbox{при}\; t\to \theta\quad\mbox{и}\quad \Big(\frac{|J_n|}{\sqrt{I_n}}\Big)^{\alpha/(1+\alpha)}\big|\overline {M_n'}(\theta)\big|=O_p(1).
\end{equation}

В частности, в двух крайних случаях $\alpha=0$ и $\alpha=1$ условия (\ref{p1-415!}) и (\ref{p1-416!}) преобразуются к следующему виду:
$$
\limsup\limits_{n\to\infty}\Big|\frac{J_n(t)}{J_n}-1\Big|=o(1)\quad\mbox{при}\quad t\to \theta\qquad \mbox{и}\qquad\overline {M_n'}(\theta)\stackrel{p}{\to}0, $$
либо
$$
 \limsup\limits_{n\to\infty}\Big|\frac{J_n(t)}{J_n}-1\Big|=O(|t-\theta|)\quad\mbox{при}\quad t\to \theta\qquad \mbox{и}\qquad \frac{|J_n|^{1/2}}{I_n^{1/4}}\big|\overline {M_n'}(\theta)\big|=O_p(1).
$$ \hfill$\square$
\end{rem}

{\bf 2.3. } Одношаговые $M$-оценки принято использовать в качестве явных приближений для состоятельных $M$-оценок, т.е. для состоятельных  статистик $\widetilde\theta_n,$ являющихся решением уравнения
\begin{equation} \label{p1-10}
M_n(t)=0 \qquad \mbox{при}\quad M_n(t)=\sum\limits_{i=1}^n M_i(t,X_i).
\end{equation}
Тем не менее, как уже  упоминалось  во введении,  одношаговые $M$-оценки и в случае разнораспределенных наблюдений
 сохраняют отмеченное, например, в [\ref{1956-L}], [\ref{1975-B}], [\ref{1991-L}], [\ref{2014-F}], [\ref{2014-1}], [\ref{2015-1?}] и [\ref{2222}]    в случае однородных выборочных данных  замечательное и несколько парадоксальное  свойство: одношаговые
   оценки могут обладать  нужными свойствами вне зависимости от того, обладают ли этими свойствами собственно сами $M$-оценки.

 Обозначим через $\widetilde\theta_n(t)$ ближайшее к $t$ решение уравнения (\ref{p1-10}).
 Нам также потребуется следующий вариант закона больших чисел
\begin{equation} \label{p1-11}
{\sum\limits_{i=1}^n\sup\limits_{t: |t-\theta|\leq \delta}M_i'(t,X_i)}\Big/{\sum\limits_{i=1}^n
{\bf E}\sup\limits_{t: |t-\theta|\leq \delta}M_i'(t,X_i)}\stackrel{p}{\to}1
\end{equation}
при некотором $\delta>0$.
 Имеет место следующее утверждение о поведении $\widetilde\theta_n(\theta)$.

\begin{theor}\label{p1-t3}
Пусть  выполнены условия   $(A_1)$-$(A_3)$,
   величины $J_n$  одного знака  при всех достаточно больших~$n$, при некотором  $0<\delta_\theta\leq
\infty$
\begin{equation} \label{p1-11-}
\limsup\limits_{n\to\infty}\overline\tau_n(\delta_\theta)<1
 \end{equation}
 и соотношение  {\rm(\ref{p1-11})} имеет место  при $\delta=\delta_\theta$.

Тогда с вероятностью, стремящейся к $1$, на интервале $(\theta-\delta_\theta,\theta+\delta_\theta)$ функция $M_n(t)$ из {\rm(\ref{p1-10})} строго монотонна {\rm(}убывает при $J_n<0$ и возрастает при $J_n>0${\rm)} и меняет знак, т.е.
  с вероятностью, стремящейся к $1$, на этом интервале
существует и единственно решение $\widetilde\theta_n(\theta)$ уравнения {\rm (\ref{p1-10})}.

Если дополнительно выполнено условие~$(A_4)$ и сходимость {\rm(\ref{p1-11})} имеет место  при всех $\delta=\delta_\theta$ таких, что верно {\rm(\ref{p1-11-})},   то
 \begin{equation}\label{p1-38}
\frac{J_n}{\sqrt{I_n}}\big(\widetilde\theta_{n}(\theta)-\theta\big)
\Longrightarrow {\cal N}\big(0,1\big).
 \end{equation}
\end{theor}

Таким образом, согласно теореме \ref{p1-t3},  при весьма  слабых
 ограничениях
 существует корень  $\widetilde\theta_n(\theta)$ уравнения (\ref{p1-10}),  удовлетворяющий (\ref{p1-38}). Этот корень уравнения (\ref{p1-10}) и приближают одношаговые $M$-оценки.

  Стоит отметить, что  $\widetilde\theta_n(\theta)$ --- это не статистика, поскольку эта величина  зависит не только от выборки, но и от неизвестного параметра $\theta$.
 Однако в случае, когда известна некоторая
 предварительная состоятельная  оценка $\theta_n^*$, то суперпозиция $\widetilde\theta_n(\theta_n^*)$ (которая является статистикой) будет обладать
  свойством:    ${\bf P}(\widetilde\theta_n(\theta)=\widetilde\theta_n(\theta_n^*))\to 1$ с ростом~$n$.
  Иными словами, алгоритм построения оценки состоит в следующем: мы выбираем ближайший к состоятельной оценке $\theta_n^*$ корень $M$-уравнения.
  Далее, если  вычисление    $\widetilde\theta_n(\theta_n^*)$
 затруднительно, то процедура одношагового $M$-оценивания --- удобная альтернатива,  позволяющая  в явном виде определять приближения для
 $\widetilde\theta_n(\theta_n^*)$.

Отметим еще, что если существует состоятельная  $M$-оценка
  $\widetilde\theta_n$, то
 \begin{equation}\label{p1-40}
  {\bf P}\big(\widetilde\theta_{n}=\widetilde\theta_n(\theta)\big)\to
  1.
 \end{equation}
Тем самым из теоремы \ref{p1-t3} и сходимости (\ref{p1-40})   можно извлечь
 следующее утверждение об асимптотическом поведении  собственно $M$-оценок в случае разнораспределенных наблюдений.

\begin{cor} \label{p1-an}
 Пусть выполнены условия $(A_1)$-$(A_4)$, величины $J_n$ одного знака при всех достаточно больших $n$,
 сходимость  {\rm(\ref{p1-11})} имеет место для всех $\delta=\delta_\theta $ таких,
  что справедливо соотношение {\rm(\ref{p1-11-})},
  и $M$-оценка $\widetilde \theta_n$ состоятельна.
Тогда
 \begin{equation}\label{p1-za7}
\frac{J_n}{\sqrt{I_n}}\big(\widetilde \theta_{n}-\theta\big)\Longrightarrow {\cal
N}\big(0,1\big).
\end{equation}
\end{cor}

\begin{rem} \label{p1-zam7}
В [\ref{1979-Y}] применительно к робастному оцениванию в линейных моделях   отмечается следующее интересное свойство одношаговых $M$-оценок: если
в качестве предварительной оценки параметра выбрать собственно саму $M$-оценку, то одношаговая $M$-оценка совпадает с $M$-оценкой.
Другими словами, если $\theta_n^*\equiv \widetilde\theta_n,$ то $\theta_{n,M}^{**}\equiv \widetilde\theta_n.$
Используя этот эффект, в  [\ref{1979-Y}] асимптотическая нормальность $M$-оценок специального вида получена из утверждений об  асимптотической нормальности одношаговых $M$-оценок  в [\ref{1975-B}].

Здесь данный прием позволяет из  теоремы \ref{p1-t1} получить следующее утверждение: если выполнены условия $(A_1)$-$(A_4)$ и $\frac{J_n}{\sqrt{I_n}}(\widetilde\theta_n-\theta)=O_p(1),$
то имеет место сходимость (\ref{p1-za7}).

Отметим, что в следствии \ref{p1-an},  извлеченном  из теоремы о поведении $\widetilde\theta_n(\theta)$, для справедливости  (\ref{p1-za7}) относительно точности $M$-оценки $\widetilde \theta_n$ достаточно лишь  условия состоятельности. \hfill$\square$
\end{rem}

\begin{rem} \label{p1-z7}
При доказательстве асимптотической нормальности $M$-оценок
 в [\ref{2007-B},  \S65, теорема 3] среди прочих технических ограничений предполагается
 (в обозначениях настоящей работы), что
 функции $M_i(t,x)$ непрерывно дифференцируемы по первому аргументу,
$|J_n|/\sum\limits_{i=1}^n{\bf
E}\sup\limits_{t\in\Theta}|M'_i(t,X_i)|>\Delta>0 $ и
\begin{equation*}
\frac{\sum\limits_{i=1}^n{\bf
E}\sup\limits_{t\in\Theta,t+u\in\Theta,|u|\leq\delta}|M'_i(t+u,X_i)-M_i'(t,X_i)|}
{\sum\limits_{i=1}^n{\bf
E}\sup\limits_{t\in\Theta}|M'_i(t,X_i)|}\leq \omega(\delta)\to
0\quad\mbox {при}\quad\delta\to 0,
  \end{equation*}
где $\omega(\cdot)$ не зависит от $n$. Таким образом,
условие $(A_4)$ следствия~\ref{p1-an} оказывается несколько слабее,
нежели соответствующие ограничения из [\ref{2007-B}]. Кроме того, в
 [\ref{2007-B}] предполагается, что $\theta$ является единственным решением уравнения
 $\sum\limits_{i=1}^n {\bf E}M_i(t,X_i)=0$ (здесь распределение $X_i$ зависит от $t$).
\hfill$\square$
\end{rem}

\begin{rem} \label{p1-z10}%
  В данной работе  мы сосредоточили свое  внимание на асимптотической нормальности одношаговых $M$-оценок.  Направление же исследований, связанное со скоростью сближения одношаговых $M$-оценок и собственно состоятельных $M$-оценок,
  лежит за рамками данного исследования.
   Тем не менее, отметим, что, на наш взгляд,   вместо изучения скорости сближения $\theta_{n,M}^{**}$ и
 $\widetilde\theta_n$ (см. (\ref{p0-4}) и  уточнения далее)
  может быть оправдано  (как в случае одинаково распределенной, так и разнораспределенной выборки)
   доказывать соотношения вида $\sqrt{n}\big(\theta_{n,M}^{**}-\widetilde\theta_n(\theta)\big)\stackrel{p}{\to}0$  и
$\theta_{n,M}^{**}-\widetilde\theta_n(\theta)=O_p(n^{-2\tau})$ при $\tau\in(1/4,1/2]$
(либо заменить $\widetilde\theta_n(\theta)$ на
$\widetilde\theta_n(\theta_n^*)$). Иными словами, может быть полезным
 изучать сближение одношаговой $M$-оценки и  ближайшего к $\theta$ корня $\widetilde\theta_n(\theta)$ уравнения (\ref{p1-10}) (или ближайшего корня этого уравнения к некоторой состоятельной оценке $\theta_n^*$), вместо  сближения одношаговой $M$-оценки и состоятельной $M$-оценки $\widetilde\theta_n$.
Такой подход позволяет  расширить  область практического применения одношаговых $M$-оценок.\hfill$\square$
\end{rem}

В завершении раздела приведем два примера.
\begin{ex}\label{p1-prim!}{\rm
В этом примере   два корня $M$-уравнения определены с вероятностью, стремящейся к $1$, при этом  состоятельной $M$-оценки не существует, но
ближайший к $\theta$ корень $\widetilde \theta_n(\theta)$ почти наверное  сходится  к $\theta$.
Итак, пусть  $X_1,\ldots,X_n$ --- выборка   одинаково распределенных наблюдений из параметрического семейства
${\cal L}_\theta(\cdot)=\theta L^{(1)}(\cdot)+\theta^2L^{(2)}(\cdot)+(1-\theta-\theta^2)L^{(3)}(\cdot),$ где распределения $L^{(k)},$ $k=1,2,3$,  не известны, но известны первые  моменты    и
 $\Theta=\{\theta>0: \theta^2+\theta<1\}\equiv \{\theta: 0<\theta<(\sqrt{5}-1)/2\}$, $\theta\neq 1/2$. Положим
$$M(\theta,X_i)=X_i-\theta a_1-\theta^2 a_2-(1-\theta-\theta^2)a_3,\quad \mbox{где }\quad a_k=\int\limits_{\mathbb{R}}tdL^{(k)}(t).$$
Соответствующее  $M$-уравнение имеет вид
$$0=\sum\limits_{i=1}^n M(\theta,X_i)=n\overline X-n\big(\theta a_1+\theta^2 a_2+(1-\theta-\theta^2)a_3\big),\quad \overline X= \sum\limits_{i=1}^nX_i/n.$$
Корни  $\widetilde\theta_{n1}<\widetilde\theta_{n2}$  квадратного уравнения
$\theta^2(a_2-a_3)+\theta(a_1-a_3)+a_3-\overline X=0$
есть
$$
\widetilde\theta_{n1,2}=\frac{a_3-a_1\pm\sqrt{(a_1-a_3)^2-4(a_2-a_3)(a_3-\overline X)}}{2(a_2-a_3)}.
$$
Выберем  распределения $L^{(k)}$  так, чтобы $a_3-a_1=a_2-a_3=\delta>0.$ В этом случае с вероятностью $1$
$$
a_3-\overline X\to (\theta-\theta^2)\delta,\quad
\widetilde\theta_{n1,2}=\frac{1}{2}\pm \frac{1}{2}\sqrt{1-4(a_3-\overline X)/\delta}{\to}\frac{1}{2}\pm \frac{1}{2}|1-2\theta|.
$$
 Асимптотически дискриминант квадратного  уравнения равен $(1-2\theta)^2> 0$, т.е.  корни $\widetilde\theta_{n1}$ и $\widetilde\theta_{n2}$ определены с вероятностью, стремящейся к $1$ и  согласно вышеприведенному определению  являются $M$-оценками. При этом эти $M$-оценки
 не являются состоятельными оценками параметра $\theta$, поскольку почти наверное  $\widetilde\theta_{n1}\to\theta$ лишь при $0<\theta<1/2$ и
 $\widetilde\theta_{n2}\to\theta$, если только $1/2< \theta\leq(\sqrt{5}-1)/2$.
  Ближайший к $\theta$ корень $\widetilde \theta_n(\theta)$    при достаточно больших $n$ определяется  равенством
 $$
 \widetilde \theta_n(\theta)=\begin{cases}
\widetilde\theta_{n1},\quad \mbox{если }\; 0<\theta< 1/2,
\\
\widetilde\theta_{n2},\quad \mbox{если}\; 1/2< \theta\leq(\sqrt{5}-1)/2
\end{cases}
 $$
 и сходится с вероятностью $1$ к $\theta$.  \hfill$\square$

}
\end{ex}

\begin{ex}\label{p1-pr-L}{\rm
Примеры параметрических семейств распределений, для которых   оценка максимального правдоподобия   не существует или не состоятельна, а оценки Фишера   в известном смысле асимптотически эффективны, в случае одномерного параметра  построены в работах Л.Ле Кама  и достаточно сложны.
В случае многомерного  параметра подобные примеры  строить  достаточно просто. Например, равновесная смесь стандартного нормального распределения и
нормального распределения с неизвестными параметрами $\alpha$ и $\sigma^2$, т.е.  плотность распределения есть
 $$
f(\tha,x)=\frac{1}{2}\varphi\big((x-\alpha)/\sigma\big)+\frac{1}{2}\varphi(x),\qquad \mbox{где }\quad \varphi(x)=e^{-x^2/2}/\sqrt{2\pi},\quad\tha=(\alpha,\sigma^2).
 $$ С~одной стороны,  оценка максимального правдоподобия  параметра $\tha=(\alpha,\sigma^2)$ здесь не существует,
 поскольку
 функция правдоподобия неограниченно возрастает в окрестностях $n$ двумерных точек $(X_i,0)$, $i=1,\ldots,n$,  лежащих на границе  открытой полуплоскости $\{\alpha\in \mathbb R,\, \sigma>0\}$, представляющей собой параметрическое множество. С другой стороны,  соответствующие оценки Фишера в этой задаче асимптотически эффективны!

Аналогичная ситуация имеет место  и в более общей ситуации, когда  плотность распределения выборочных данных задается соотношением
$$
f(\tha,x)=\frac{p}{\sigma}g\big((x-a)/\sigma\big)+\frac{1-p}{\tau}g\big((x-b)/\tau\big),\qquad g(0)\neq 0, \quad \tha=(p,a,b,\sigma,\tau).
$$
В качестве  $\theta_n^*$ здесь можно выбрать оценку по методу моментов.
Данный пример содержится, например,  в  монографии Э. Лемана [\ref{1991-L}].
%
 \hfill$\square$

}
\end{ex}

\section{Приложения к задачам регрессии }

{\bf 3.1.} В этом разделе  приведем несколько относящихся к регрессионному анализу  примеров возможного использования полученных  результатов.
Рассмотрим классическую  модель нелинейной регрессии. В этом случае параметр $\theta$ нужно оценить по наблюдениям $X_1,\ldots,X_n,$ имеющим следующую структуру:
\begin{equation}\label{p1-601-}
X_i=f_i(\theta)+\varepsilon_i, \qquad {\bf E}\varepsilon_i=0,\qquad {\bf D}\varepsilon_i=\sigma^2,\qquad i=1,\ldots,n,
\end{equation}
где $\{f_i(\cdot)\}$  --- известные функции,  погрешности $\{\varepsilon_i\}$ --- последовательность независимых
 случайных величин,  параметр  $\sigma^2$
может быть не известным.

Оценка метода наименьших квадратов (МНК) при некоторых условиях регулярности определяется здесь  как решение $\widetilde\theta_n$  уравнения
\begin{equation}  \label{p1-226}
\sum\limits_{i=1}^n M_i\big(t,X_i\big)=0 \quad \mbox{при}\quad M_i(t,X_i)=f'_i(t)(X_i-f_i(t)),
 \end{equation}
 минимизирующее сумму квадратов погрешностей:  $\widetilde\theta_n:=\arg\min\limits_\theta\sum\limits_{i=1}^n(X_i-f_i(\theta))^2.$
При широких ограничениях МНК-оценка $\widetilde\theta_n$ оказывается асимптотически нормальной с асимптотической дисперсией $\sigma^2\big/\sum\limits_{i=1}^n\big(f_i'(\theta)\big)^2$.
 Обоснованием точности нелинейного МНК в случае, когда наблюдения имеют нормальное распределение, может служить совпадение в этом случае МНК-оценок  и оценок метода  максимального правдоподобия. Трудности, связанные  с  вычислением МНК-оценок в нелинейной регрессии,   зачастую
   обусловлены наличием  большого числа локальных минимумов у минимизируемой  функции, а потому
 при неудачном выборе начального приближения параметра итерационные процедуры  обнаруживают лишь локальный минимум, ближайший к этой стартовой точке (см., например, [\ref{1989-D}]).
   Методология  одношагового улучшения оценок в случаях,  когда возможно отыскать некоторую предварительную оценку для $\theta$, позволяет в явном виде строить приближения для решений уравнения (\ref{p1-226}).

\begin{ex}\label{p1-e2}
{\em
Пусть регрессионная модель задается соотношениями:
\begin{equation}  \label{p1-233}
X_i=\sqrt{1+z_i\theta}+\varepsilon_i, \quad {\bf E}\varepsilon_i=0,\quad {\bf D}\varepsilon_i=\sigma^2,\quad i=1,\ldots,n,
\end{equation}
где числа $\{z_i>0\}$ известны.
 В монографии [\ref{1989-D}] указано, вычисление   МНК-оценки обычными  итерационными методами для этой модели (а также для регрессионной функции из приводимого далее примера \ref{p1-en-}) весьма трудоемко.
  В   [\ref{2014-S}] в качестве альтернативы для   оптимальной, но трудновычислимой  МНК-оценки   построена следующая асимптотически нормальная, но не асимптотически оптимальная (в смысле точности МНК) оценка для параметра $\theta>0$:
\begin{equation}  \label{p1-234}
\theta_n^*=\sum\limits_{i=1}^n c_{i}X_i^2\Big/\sum\limits_{i=1}^n c_{i}z_i,\quad \mbox{где}\quad \sum\limits_{i=1}^n c_{i}=0.
\end{equation}

Используя конструкцию одношаговых приближений, точность оценивания из [\ref{2014-S}]  можно  улучшить, при этом  этом статистика (\ref{p1-234}) может  использоваться в качестве предварительной оценки.
Для функций из (\ref{p1-226}), имеем
\begin{equation}\label{p1-234+}
\begin{split}
M_i'(\theta,X_i)=f_i''(\theta)X_i-\big(f_i'(\theta)f_i(\theta)\big)_\theta',\\ {\bf E}M_i'(\theta,X_i)=-(f_i'(\theta))^2
=f_i''(\theta)f_i(\theta)-\big(f_i'(\theta)f_i(\theta)\big)_\theta'.
\end{split}
\end{equation}
Поскольку для данной модели $f_i'(\theta)f_i(\theta)$ не зависит от $\theta$, то используя определения (\ref{p0-0}), (\ref{p1-226})   и первое равенство в (\ref{p1-234+}), получаем
\begin{gather*}
\theta_{n,M}^{**}=\theta_n^*+{2\sum\limits_{i=1}^n\frac{z_i}{\sqrt{1+z_i\theta}}\big(X_i-\sqrt{1+z_i\theta}\big)}\Big/{\sum\limits_{i=1}^n
 z_i^2(1+z_i\theta_n^*)^{-3/2}X_i}.
 \end{gather*}
 Оценки $\theta_{n,J}^{**}$ из (\ref{p1-777}) определяются  не однозначно: все зависит от выбора оценки $J_n^*$ для $J_n\equiv\sum\limits_{i=1}^n{\bf E}M_i'(\theta,X_i).$ Например, следуя второму равенству в  (\ref{p1-234+}), здесь  можно положить
 $J_n^*=-\sum\limits_{i=1}^n (f_i'(\theta_n^*))^2$ или
$J_n^*=\sum\limits_{i=1}^n f_i''(\theta_n^*)\big(rf_i(\theta_n^*)-(r-1)X_i\big)$ или $J_n^*=\sum\limits_{i=1}^n f_i''(\theta_n^*)\big(rX_i -(r-1)f_i(\theta_n^*)\big)$ при любом $r\in \mathbb{R}$. В первом случае одношаговая  оценка $\theta_{n,J}^{**}$ примет следующий вид:
\begin{gather*}
\theta_{n,J}^{**}=\theta_n^*+{2\sum\limits_{i=1}^n\frac{z_i}{\sqrt{1+z_i\theta}}\big(X_i-\sqrt{1+z_i\theta}\big)}\Big/{\sum\limits_{i=1}^n
 z_i^2(1+z_i\theta_n^*)^{-1}}.
 \end{gather*}
 Асимптотическая дисперсия всех предложенных  одношаговых $M$-оценок совпадает с асимптотической дисперсией МНК-оценки. Эта оптимальная асимптотическая дисперсия есть
$I_n/J_n^2=\sigma^2\big/\sum\limits_{i=1}^n\big(f_i'(\theta)\big)^2\equiv 4\sigma^2\Big/\sum\limits_{i=1}^nz_i^2(1+\theta z_i)^{-1}$.
~\hfill$\square$
}
\end{ex}

 Рассмотрим теперь гетероскедастическую модель нелинейной регрессии,   в которой дисперсии наблюдений зависят не только от номера наблюдения, но и от основного неизвестного параметра. В этом случае параметр $\theta$ нужно оценить по наблюдениям $X_1,\ldots,X_n,$ имеющим следующую структуру:
\begin{equation}\label{p1-601-!}
X_i=f_i(\theta)+\varepsilon_i, \qquad {\bf E}\varepsilon_i=0,\qquad {\bf D}\varepsilon_i=\sigma^2/w_i(\theta),\qquad i=1,\ldots,n,
\end{equation}
где $\{f_i(\cdot)\}$ и $\{w_i(\cdot)\}$ --- известные функции.
Оптимальная
в некотором классе $M$-оценка $\widetilde \theta_n$ (так называемая оценка квази-правдоподобия; см., например, [\ref{1997-H}])  определяется здесь как решение уравнения
\begin{equation}\label{p1-602}
\sum\limits_{i=1}^nw_i(t)f_i'(t)(X_i-f_i(t))=0.
\end{equation}
 При широких ограничениях оценка квази-правдоподобия  $\widetilde\theta_n$ оказывается асимптотически нормальной с асимптотической дисперсией $\sigma^2\big/\sum\limits_{i=1}^n w_i(\theta)\big(f_i'(\theta)\big)^2$.
Подробнее одношаговые процедуры для  нахождения в явном виде оценок, имеющих асимптотически точность оценок квази-правдоподобия,  рассмотрены автором в [\ref{2015-3?}].

Предположим, дополнительно к (\ref{p1-601-!}),  что наблюдения $\{X_i\}$ --- нормально распределенные величины. В этом случае оценка
 максимального правдоподобия оказывается несколько точнее, нежели оценка квази-правдоподобия. Обычно параметр $\sigma^2$ предполагается {\it не известным}. Здесь для простоты рассмотрим модельную задачу, считая,  что значение $\sigma^2$ {\it известно}. В противном случае принципиально все обсуждаемые здесь эффекты сохранятся, но нужно воспользоваться уже   двумерными  аналогами приведенных в  \S2  результатов (см. \S4).

Для модели (\ref{p1-601-!}) оценка максимального правдоподобия $\widetilde\theta_n$ является  решением следующего уравнения:
\begin{equation}\label{p1-599-}
\begin{split}
\sum\limits_{i=1}^n M_i(t,X_i)=0, \quad \\ M_i(t,X_i)= \frac{(X_i-f_i(t))f_i'(t)w_i(t)}{\sigma^2}-
\frac{(X_i-f_i(t))^2w_i'(t)}{2\sigma^2}+\frac{w_i'(t)}{2w_i(t)}.
\end{split}
\end{equation}
Одношаговая оценка $\theta_{n,J}^{**}$, определяемая соотношением (\ref{p1-777}) c функциями $\{M_i(\cdot,\cdot)\}$ из (\ref{p1-599-}),
примет следующий вид:
\begin{equation}\label{p1-578-}
\theta_{n,J}^{**}=\theta_n^*+\frac{\displaystyle\sum\limits_{i=1}^n\left((X_i-f_i(\theta_n^*))f_i'(\theta_n^*)w_i(\theta_n^*)-
\frac{(X_i-f_i(\theta_n^*))^2w_i'(\theta_n^*)}{2}+\frac{\sigma^2w_i'(\theta_n^*)}{2w_i(\theta_n^*)}\right)}{\displaystyle\sum\limits_{i=1}^n \displaystyle\left(\big(f_i'(\theta_n^*)\big)^2w_i(\theta_n^*)+\sigma^22^{-1}\big(w_i'(\theta_n^*)/w_i(\theta_n^*)\big)^2\right)}.
\end{equation}
Отметим, что можно здесь воспользоваться и статистикой $\theta_{n,M}^{**}$ из (\ref{p1-7}), но представление для $M_i'(t,X_i)$ оказывается  более громоздким, нежели информация Фишера, соответствующая наблюдению $X_i$ и участвующая в знаменателе представления (\ref{p1-578-}) для оценки $\theta_{n,J}^{**}$.

 Таким образом, одношаговые $M$-оценки позволяют в явном виде находить оценки, имеющие асимптотически ту же точность, что и трудно вычислимая в данной ситуации оценка максимального правдоподобия  $\widetilde\theta_n$. В частности, в данной задаче при выполнении условий теоремы \ref{p1-t1!-!!} оценка
 $\theta_{n,J}^{**}$, определенная в (\ref{p1-578-}), является асимптотически нормальной с асимптотической дисперсией
\begin{equation}\label{p1-598-}
I_n/J_n^2\equiv \sigma^2\Big/\sum\limits_{i=1}^n\left( w_i(\theta)\big(f_i'(\theta)\big)^2+\frac{\sigma^2}{2}(w_i'(\theta)/w_i(\theta))^2\right).
\end{equation}

\begin{ex} \label{p1-e1-}
{\em
Рассмотрим одномерную
 модель нелинейной регрессии, порожденную известным биохимическим  уравнением  Михаэлиса--Ментен (см., например, [\ref{2011-M}], [\ref{2009-D}], [\ref{2001-2}]). В этом случае  наблюдения $\{X_i\}$
представимы в виде
\begin{equation*} 
X_i=\frac{a_i}{1+b_i\theta}+\varepsilon_i,
\qquad   i=1,\ldots,n,
 \end{equation*}
где $\{a_i>0\}$ и $\{b_i>0\}$ --- известные числовые последовательности.
 И пусть  выполнено (\ref{p1-601-!}), а наблюдения имеют нормальные распределения.
В качестве оценки $\theta_n^*$ для построения $\theta_{n,J}^{**}$ можно использовать следующую статистику (асимптотически нормальную при широких ограничениях), построенную и исследованную в [\ref{2000-1}]:
\begin{equation}  \label{p1-227g+}
\theta_n^*=\sum\limits_{i=1}^nc_i(a_i-X_i)\Big/\sum\limits_{i=1}^nc_ib_iX_i,
 \end{equation}
 где $\{c_i\}$ --- некоторые числа. Тот факт, что асимптотическая дисперсия оценки $\theta_n^*$ из (\ref{p1-227g+}), равная
  $ \sigma^2\sum\limits_{i=1}^n c_i^2(1+b_i\theta)^2w_i^{-1}(\theta)\Big/\Big( \sum\limits_{i=1}^n c_ia_ib_i(1+b_i\theta)^{-1}\Big)^2,$ при любом выборе констант $\{c_i\}$ больше, чем $\sigma^2\Big/\sum\limits_{i=1}^n w_i(\theta)a_i^2b_i^2/(1+b_i\theta)^2$,  установлен в [\ref{2000-1}].
 Таким образом, асимптотическая дисперсия одношаговой оценки $\theta_{n,J}^{**}$ из (\ref{p1-578-}), равная
\begin{equation*}
I_n/J_n^2\equiv \sigma^2\Big/\sum\limits_{i=1}^n\left( w_i(\theta)a_i^2b_i^2/(1+b_i\theta)^2+\frac{\sigma^2}{2}(w_i'(\theta)/w_i(\theta))^2\right),
\end{equation*}
оказывается меньше, нежели асимптотическая дисперсия предварительной  оценки $\theta_n^*$.
\hfill$\square$.
}
\end{ex}

\begin{ex}
\label{p1-en-}
{\em
Пусть наблюдения $X_i$ имеют следующую структуру:
\begin{equation}  \label{p1-230-}
X_i=a_i\theta+b_ig(\theta)+\varepsilon_i,
\qquad i=1,\ldots,n,
\end{equation}
где $\{a_i\}$ и $\{b_i\}$ --- известные числовые последовательности.
  В  [\ref{2013-E}]   построена следующая явная статистика для параметра $\theta$:
\begin{equation}  \label{p1-231-}
\theta_n^*=\sum\limits_{i=1}^n c_{i}X_i\Big/\sum\limits_{i=1}^n c_{i}a_i,
\end{equation}
где числа  $\{c_{i}\}$ таковы, что $\sum\limits_{i=1}^n c_{i}b_i=0$ и $\sum\limits_{i=1}^n c_{i}a_i\neq0.$ В
 [\ref{2013-E}]  доказано, что при широких  ограничениях оценка $\theta_n^*$ является асимптотически нормальной с асимптотической дисперсией $ \sigma^2\sum\limits_{i=1}^n c_{i}^2{\bf D}\varepsilon_i
 \Big/\Big( \sum\limits_{i=1}^n c_{i}a_i\Big)^2$.

Используя обсуждаемый здесь подход одношагового улучшения точности оценивания, оценку $\theta_n^*$ из  [\ref{2013-E}] можно улучшить.
Аналогично примеру \ref{p1-e2}, в случае однородных погрешностей $\{\varepsilon_i\}$ (см. (\ref{p1-601-})) можно использовать различные одношаговые оценки, имеющие асимптотическую точность  МНК-оценки. Если же выполнено [\ref{p1-601-!}] --- можно находить одношаговые улучшения, имеющие точность  оценок квази-правдоподобия.

В случае (\ref{p1-601-!}) и нормально распределенных наблюдений в качестве улучшения для $\theta_n^*$ из  [\ref{2013-E}] следует выбрать   одношаговую  оценку
$\theta_{n,J}^{**}$, определенную в  (\ref{p1-598-}).
В силу (\ref{p1-598-}) асимптотическая дисперсия этой   одношаговой  $M$-оценки
  определяется равенством  $I_{n}/J_{n}^2\equiv\sigma^2\big/\sum\limits_{i=1}^n\left(w_i(\theta)\big(a_i+b_ig'(\theta)\big)^2+
  \frac{\sigma^2}{2}\big(w_i'(\theta)/w_i(\theta)\big)^2
  \right).
   $
    Нетрудно видеть, что ни при каком выборе констант $\{c_{i}\}$ в (\ref{p1-231-}) для оценки $\theta_n^*$ не может быть достигнута эта точность.
        \hfill$\square$
}
\end{ex}

{\bf 3.2.} 
 %
%
%
  Обсудим вопрос о   некоторых    способах построения предварительных оценок в задачах нелинейной регрессии.
  Среди многообразия моделей нелинейной регрессии выделим, прежде всего, так называемые  {\it внутренне линейные} (или {\it внешне нелинейные}) модели регрессии.   Напомним (см., например, [\ref{1987-D2}]), что  к подобным   моделям принято относить модели нелинейной регрессии, которые теми или иными преобразованиями   можно привести к линейным моделям регрессии, содержащим изначально оцениваемые параметры.

Так, например, к этому классу относятся    модели  регрессии  (\ref{p1-233}) и  (\ref{p1-230-}). Тем  самым,  можно предложить и   другой способ  построения явных  оценок в этих   моделях,  
нежели приемы, используемые  в [\ref{2013-E}],   [\ref{2014-S}], а также в  [\ref{2015-K}].

   Действительно, для регрессионной функции из (\ref{p1-230-}) достаточно положить
        $$
    \widetilde X_i=\frac{X_{2i}}{b_{2i}}-\frac{X_{2i-1}}{b_{2i-1}}, \quad \widetilde a_i=\frac{a_{2i}}{b_{2i}}-\frac{a_{2i-1}}{b_{2i-1}}, \quad
    \widetilde\varepsilon_i=\frac{\varepsilon_{2i}}{b_{2i}}-\frac{\varepsilon_{2i-1}}{b_{2i-1}},\quad i=1,\ldots,n/2
    $$
    (без ограничения общности можно считать, что  $b_j\neq 0$ при всех $j$)
    и в модели линейной регрессии
    $\widetilde X_i=\widetilde a_i\theta+\widetilde\varepsilon_i,$ $i=1,\ldots,n/2$,  оценить $\theta$ взвешенным МНК с некоторыми весами.

Для модели (\ref{p1-233}) имеем    $X_i^2= 1+\theta z_i +\varepsilon_i^2+2\sqrt{1+\theta z_i}\varepsilon_i$. Следовательно,  параметр $\theta$ можно  оценить исходя  из  модели линейной регрессии
     $     \widetilde X_i=\theta_0+\theta  z_i+\widetilde \varepsilon_i$,
     где $\widetilde X_i=X_i^2-1,$ $\quad \widetilde \varepsilon_i=2\sqrt{1+\theta z_i}\varepsilon_i+\varepsilon_i^2-{\bf E}\varepsilon_1^2,$ $\quad\theta_0={\bf E}\varepsilon_1^2$. В обоих примерах   ${\bf E}\widetilde\varepsilon_i=0$  и
  при широких ограничениях  оценки взвешенного МНК c некоторыми весами
  окажутся асимптотически нормальными (но не оптимальными) для параметра $\theta$ в полученных  гетероскедастических моделях линейной регрессии.
Аналогичным образом можно оценить и  многомерный параметр $\theta_1,\ldots,\theta_m$ для более общей внутренне линейной модели $X_i=\sqrt{\theta_{1}z_{1i}+\ldots+\theta_m z_{mi}}+\varepsilon_i$ с соответствующими  моментными ограничениями на погрешности.

В качестве еще одного примера внутренне линейной модели рассмотрим случай, когда наблюдения $\{X_i\}$ имеют структуру
$X_i=\ln (1+\theta z_i)+\varepsilon_i$ с условием ${\bf E}e^{\varepsilon_i}\equiv\sigma^2<\infty$ при всех $i$.
Преобразуем это регрессионное уравнение  следующим образом: $e^{X_i}=(1+\theta z_i)\big(e^{\varepsilon_i}-{\bf E}e^{\varepsilon_1}\big)+(1+\theta z_i)
{\bf E}e^{\varepsilon_1}.$   Это есть  модель линейной регрессии $\tilde X_{i}=\theta_0+\theta_1  z_{i}+\tilde\varepsilon_{i}$, если  положить
$  \tilde X_{i}=e^{X_i}$, $\theta_0={\bf E}e^{\varepsilon_1},$ $\theta_1=\theta\theta_0$, $ \tilde \varepsilon_{i}=(1+\theta z_i)\big(e^{\varepsilon_i}-{\bf E}e^{\varepsilon_1}\big)$.
 Пусть теперь   $\theta_0^*$ и $\theta_1^*$ --- оценки взвешенного МНК, а искомую оценку параметра $\theta$ определим соотношением   $\theta^*= \theta_1^*/\theta_0^*$.  Нетрудно проверить, что при широких ограничениях $\theta_n^*$  является асимптотически нормальной оценкой.
Описанный выше прием очевидным образом обобщается на случай многомерного параметра, когда отклики имеют структуру
$X_i=\ln(1+\theta_{1}z_{1i}+\ldots+\theta_m z_{mi})+\varepsilon_i.$
 Несколько иная  аргументация при  построении явных асимптотически нормальных оценок параметра   $\theta$ для такой логарифмической модели была предложена  в  [\ref{2015-K}].

  К сожалению, класс внутренне линейных моделей достаточно узок, и подобное свойство возможного преобразования нелинейной модели к линейной является скорее исключением, нежели правилом. Поэтому возникает необходимость разработки некоторых подходов получения предварительных оценок для собственно нелинейных моделей.

Рассмотрим  типичную регрессионную модель, когда  наблюдения $\{X_i\}$ имеют структуру
\begin{equation}\label{p1-m7}
 X_{i}=f(\theta,z_i)+\varepsilon_{i},\qquad {\bf E}\varepsilon_i=0,\quad {\bf E}\varepsilon_i^2<\infty,\quad i=1,\ldots,n,
\end{equation}
 где функция  $f(\cdot,\cdot)$ и числовая последовательность $\{z_i\}$ известны,  погрешности $\{\varepsilon_i\}$  независимы и $\theta\in \Theta=(a,b)$.

Предлагаемый подход основан на
использовании интегральных сумм Римана, а потому удобнее при каждом фиксированном $n\geq 1$ упорядочить $z_1,\ldots,z_n$ по возрастанию:
$z_{n:1}\leq \ldots\leq z_{n:n}.$ В этом случае соответствующим образом перенумеруются  $X_i$ и  $\varepsilon_i$ (с сохранением условия независимости) и модель (\ref{p1-m7}) примет  вид
\begin{equation*} 
 X_{ni}=f(\theta,z_{n:i})+\varepsilon_{ni},\qquad {\bf E}\varepsilon_{ni}=0,\quad {\bf E}\varepsilon_{ni}^2<\infty,\quad i=1,\ldots,n.
\end{equation*}

В приводимой  ниже теореме \ref{p1-co1} предполагается, что
 с ростом $n$
 расстояния между точками $z_{n:i}$ и $z_{n:i-1}$ с определенной скоростью сходятся к нулю, при этом  возможна  как ситуация, когда  все $z_i$ лежат на конечном отрезке $[c,d]$
(в этом случае совокупность  $\{z_{n:i}\}$ образует измельчающееся разбиение отрезка
 $[c,d]$), так и случай неограниченного возрастания  $z_{n:n}$.

Положим  $\Delta z_{ni}:=z_{n:i}-z_{n:i-1},\quad z_{n:0}= c,$
\begin{eqnarray}\label{p1-m88}
\omega_{f,\theta}(\Delta)=\sup\limits_{|t_1-t_2|\leq \Delta }|f(\theta,t_1)-f(\theta,t_2)|, \qquad T(t)=\displaystyle\int\limits_c^df(t,z)dz.
\end{eqnarray}

\begin{theor} \label{p1-co1}
Пусть  $z_i\in [c,d]$ при всех $i$,  $z_{n,n}\to d$ {\rm(}допускается, что~$d=~\infty${\rm)}, $\max_{i\leq n}\Delta z_{ni}\to 0$,  функция
 $f(\theta,z)$ интегрируема  по Риману
 второму аргументу на $[c,d]$,    функция $T(t)$ из~{\rm (\ref{p1-m88})}  строго монотонна и непрерывна на $(a,b)$, а обратная функция $T^{-1}(t)$ удовлетворяет условию Гельдера с показателем $p\in(0,1]$.   Кроме того, пусть
\begin{eqnarray}\label{p1-m9-}
\alpha_n^{1/p}d_n\to 0\qquad\mbox{при}\quad d_n^2=\sum\limits_{i=1}^n\big(\Delta z_{ni}\big)^2{\bf E}\varepsilon_{ni}^2,
\\\label{p1-m9}
\alpha_n^{1/p}\sum\limits_{i=1}^n\omega_{f,\theta}(\Delta z_{ni})\Delta z_{ni}\to 0,\qquad \alpha_n^{1/p}\int\limits_{z_{n:n}}^d f(\theta,z)dz\to 0.
\end{eqnarray}
 Тогда  оценка
\begin{equation*} 
\theta_n^*=T^{-1}\Big(\sum\limits_{i=1}^n\Delta z_{ni}X_{ni}\Big)
\end{equation*}
является $\alpha_n$-состоятельной.

Если дополнительно функция $T(t)$ непрерывно дифференцируема,  условия из  {\rm(\ref{p1-m9})} выполнены при замене  $\alpha_n^{1/p}$ на $1/d_n$, а вместо  условия   {\rm(\ref{p1-m9-})}
имеют место сходимости $d_n\to 0$ и  $\displaystyle\sum\limits_{i=1}^n\Delta z_{ni} \varepsilon_{ni}/d_n\Rightarrow {\cal N}(0,1)$, то
\begin{equation*} 
(\theta_n^*-\theta)/d_n \Longrightarrow{\cal N}\big(0,[T'(\theta)]^{-2}\big).
\end{equation*}
\end{theor}

Обобщая пример 2, пусть $f(\theta,z_i)= (1+\theta z_i)^r,$ $\theta>0,$ $z_i\in [0,1]$,  $r\neq 0, -1$. В этом случае
   $T(\theta)=\frac{(1+\theta)^{r+1}-1}{\theta(r+1)}$ и   значения
обратной функции $T^{-1}(\cdot)$
 могут быть легко  вычислены (с наперед заданной точностью) по таблице для функции $T(t)$.  В случае, когда  $r<-1$,
  $c=0$ и $d=\infty$ имеем
$T(\theta)=\frac{1}{\theta|r+1|}$.
Отметим, что модель $f(\theta,z_i)= (1+\theta z_i)^{1/2}$ из примера 2 является внутренне линейной только в случае, когда ${\bf E}\epsilon_i^2=\sigma^2/w_i$ и числа $\{w_i\}$ известны. Предлагаемый здесь прием позволяет находить предварительную оценку для указанной регрессионной функции   и при других (более общих) ограничениях на погрешности наблюдений.

Предлагаемый в теореме \ref{p1-co1} прием построения оценок  для функций специального вида (в частности, степенных) можно распространить
  и на случай, когда  $z_{n:n}\to \infty$,  регрессионная функция $f(\theta,z)$ не интегрируема на $[c,\infty)$, а последовательность  $\{z_{n:i}/z_{n,n}\}$ образует измельчающееся разбиение отрезка
 $[0,1]$). Кроме того, этот прием допускает обобщения на случай многомерного неизвестного параметра и многомерного регрессора. Класс таких многомерных  задач нелинейной регрессии   включает в себя, в частности,  важные в приложениях модели роста, описывающие процессы, происходящие в биологии, химии, экономике. Подробнее     методы построения предварительных $\alpha_n$-состоятельных или асимптотически нормальных (но не оптимальных) оценок  в различных  задачах  нелинейной регрессии будут рассмотрены в отдельной работе. В частности, речь идет об оценивании параметра  $\tha\in \mathbb{R}^m$ в моделях нелинейной регрессии вида $
 X_i=f(\tha,{\bf z}_i)+\varepsilon_i$,  $i=1,\ldots,n,$
 где  $f(\cdot,\cdot)$ --- некоторая известная функция, значения $k$-мерных векторов   ${\bf z}_i$ также  известны. 


{\bf 3.3.}    Иногда функции $M_i(t,X_i)$ естественным образом допускают факторизацию вида
\begin{equation} \label{p1-15w}
M_i(t,X_i)=h_i(t)\widetilde M_i(t,X_i).
\end{equation}
Например, в случае МНК-оценки (\ref{p1-226}) можно положить
$h_i(t)=f_i'(t),$ $\widetilde M_i(t,X_i)=X_i-f_i(t),$ а для 
 оценки квази-правдоподобия  (\ref{p1-602}) считать, что $h_i(t)=w_i(t)f_i'(t)$, $ \widetilde M_i(t,X_i)=X_i-f_i(t).$
  Оказывается, выделение в (\ref{p1-15w}) некоторых множителей  $h_i(\theta)$ позволяет ослаблять условия на гладкость функций, определяющих одношаговые оценки.
Ниже мы приведем простые достаточные условия для асимптотической нормальности одношаговой  $M$-оценки $\theta_{n,J}^{**}$ в случае, когда функции $h_i(t)$ и $\widetilde M_i'(t,X_i)$ удовлетворяют лишь условию Гельдера.
В частности, для двух отмеченных  примеров факторизации (\ref{p1-15w})
использование приводимой далее теоремы \ref{p1-t1!-!!w} позволяет ослабить (по сравнению с требованиями   теоремы  \ref{p1-t1!-!!})   условия  на  гладкость
 функций $f_i(t)$ и $w_i(t)$. Отметим еще, что при определении одношагового приближения $\theta_{n,J}^{**}$ для   МНК-оценки проще всего  выбрать $J_n^*=\sum\limits_{i=1}^n\big(f_i'(\theta_n^*)\big)^2$, а в случае  оценки   квази-правдоподобия  положить $J_n^*=\sum\limits_{i=1}^nw_i(\theta_n^*)\big(f_i'(\theta_n^*)\big)^2$.

Нам потребуется  следующее предположение.

$({\bf A}^-).$ Выполнено условие
$ ( A_1)$,
  при любом  $i$ на   $\Theta$  заданы   функции $h_i(t)$ и с вероятностью $1$ непрерывно дифференцируемые   по первому аргументу функции    $\widetilde M_i(t,X_i)$ такие, что для некоторых $p,q\in(0,1]$
\begin{equation} \label{p1-36w}
 \forall t_1,t_2\in\Theta\qquad |h_i(t_1)-h_i(t_2)|\leq \overline h_i|t_1-t_2|^p,
 \end{equation}
$$ \forall t\in\Theta\qquad |\widetilde M_i'(t,X_i)-\widetilde M_i'(\theta,X_i)|\leq \overline { M_i'}|t-\theta|^q, \quad {\bf E}\overline { M_i'}<\infty,
$$
при этом  ${\bf E}\widetilde M_i(\theta,X_i)=0$, ${\bf E}\widetilde M^2_i(\theta,X_i)<\infty$, ${\bf E}|\widetilde M'_i(\theta,X_i)|<\infty$,
  начиная с некоторого $n$
$
I_n=
\sum\limits_{i=1}^nh_i^2(\theta)
{\bf E}\widetilde M^2_i(\theta,X_i)>0,$ $ J_n:=
\sum\limits_{i=1}^nh_i(\theta){\bf E}\widetilde M_i'(\theta,X_i)\neq 0,
$ и имеет место сходимость $\sum\limits_{i=1}^n h_i(\theta)\widetilde M_i(\theta,X_i)/\sqrt{I_n}\Longrightarrow {\cal N}(0,1).
$
Кроме того, имеются  оценки $\theta_n^*$ и $J_n^*$, для которых 
 $ J_n^*/J_n\stackrel{p}{\to}1$ и
 \begin{eqnarray}                                                                                  \label{p1-17w}
  \begin{split}
  |\theta_n^*-\theta|^{1+p}\sum\limits_{i=1}^n\overline h_i{\bf E}\big|\widetilde M_i'(\theta,X_i)\big|\big/\sqrt{I_n}\stackrel{p}{\to} 0,\\
  |\theta_n^*-\theta|^{1+q}\sum\limits_{i=1}^n\big(|\theta_n^*-\theta|^p\overline h_i+|h_i(\theta)|{\bf E}\overline{ M_i'}\big)\big/\sqrt{I_n}\stackrel{p}{\to} 0,\qquad\\
 (\theta_n^*-\theta)(J_n^*-J_n)/\sqrt{I_n}\stackrel{p}{\to}0, \quad (\theta_n^*-\theta)\Big(\sum\limits_{i=1}^nh_i(\theta)\widetilde M_i'(\theta,X_i)-J_n\Big)\big/\sqrt{I_n}\stackrel{p}{\to}0.
  \end{split}
  \end{eqnarray}

 При выполнении условия $(A^-)$ определим одношаговую $M$-оценку соотношением
\begin{equation} \label{p1-18w}
\theta_{n,J}^{**}=\theta_n^*-\sum\limits_{i=1}^n h_i(\theta_n^*)\widetilde M_i(\theta_n^*,X_i)\big/J_n^*.
\end{equation}

\begin{theor}                                     \label{p1-t1!-!!w}
Если   выполнено  предположение   $(A^-)$ и
\begin{gather}                                                                                           \label{p1-19w}
\sum\limits_{i=1}^n\big(h_i(\theta_n^*)-h_i(\theta)\big)\widetilde M_i(\theta,X_i)\big/\sqrt{I_n}\stackrel{p}{\to}0,
 \end{gather}
 то одношаговая $M$-оценка $\theta_{n,J}^{**}$  определена с вероятностью, стремящейся к~$1$~и
\begin{gather*} 
\frac{J_n}{\sqrt{I_n}}\big(\theta_{n,J}^{**}-\theta\big)\Longrightarrow {\cal N}(0,1).
 \end{gather*}
\end{theor}

\begin{rem}
Доказательство (\ref{p1-19w}) в случае, когда функции $h_i(t)$ удовлетворяют лишь условию Гельдера,
представляет некоторую техническую трудность (см. [\ref{2015-3?}] и [\ref{2011-1}]). Без ограничения общности будем  считать, что функции $h_i(t)$
продолжены на $\mathbb{R}$ с выполнением  условия Гельдера (\ref{p1-36w}) на всей числовой прямой.
Достаточные условия для  (\ref{p1-19w}) в случае, когда  предварительная оценка $\theta_n^*$
имеет следующую структуру
\begin{equation*}\label{p1-28w}
\theta_n^*-\theta=\frac{\sum\limits_{i=1}^n u_{ni}(\theta,X_i)}{1+\sum\limits_{i=1}^n v_{ni}(\theta,X_i)}, \qquad {\bf E}u_{ni}(\theta,X_i)= {\bf E}v_{ni}(\theta,X_i)=0,\;i=1,\ldots,n,
\end{equation*}
нетрудно извлечь из теоремы 3 в  [\ref{2011-1}] (во всех примерах из раздела 3.2 оценка
 $\theta_n^*$ удовлетворяет указанному требованию из [\ref{2011-1}]).
 В частности, условие (\ref{p1-19w}) выполнено, если
$d_{nu}^2 +\sum\limits_{i=1}^n {\bf E}v_{ni}^2(\theta,X_i)\to 0$ при $d_{nu}^2:=\sum\limits_{i=1}^n {\bf E}u_{ni}^2(\theta,X_i)$ и
\begin{eqnarray}\label{p1-29w}
\begin{split}
d_{nu}^{2p}\left(\sum\limits_{i=1}^n \Big(\overline h_i^2{\bf E}\widetilde M_{i}^2(\theta,X_i)\Big)^{1/(2-p)}\right)^{2-p}\Big/I_n\to 0,
 \end{split}
 \end{eqnarray}
Отметим, что центральное условие  (\ref{p1-29w}) имеет место, если
\begin{gather}\label{p1-31w}
 n^{(1/p-1)/2}d_{nu}\to 0,\quad {1}/I_{n}=O\big(1/n\big),\quad \sup\limits_{n,i}\{\overline h_i,\; {\bf E} M_i^2(\theta,X_i)\}<\infty.
 \end{gather}
В частности, первое условие в (\ref{p1-31w}) выполнено в следующих трех случаях:
$$
\sqrt{n}d_{nu}=O(1)\; \mbox{и}\; p>1/2, \qquad n^{1/4}d_{nu}\to 0\; \mbox{и}\; p\geq 2/3,\qquad d_{nu}\to 0\; \mbox{и}\; p=1.
$$
Если  $\theta_n^*$ имеет  дробно-линейную структуру, то достаточные условия для (\ref{p1-19w}) можно получать и в терминах
сближения $\theta_n^*$ и $\theta$ по вероятности (т.е. в терминах условия  $n^{(1/p-1)/2}(\theta_n^*-\theta)\stackrel{p}{\to} 0$) (см. подробности в [\ref{2013-2}]).

В [\ref{2015-3?}] (без использования предположения о специальной структуре $\theta_n^*$) доказано, что для выполнения (\ref{p1-19w}) 
достаточно, чтобы имели место сходимости
 \begin{gather}
 \big({\bf E}|\theta_n^*-\theta|^2\big)^p\sum\limits_{i=1}^n\overline h_i^2{\bf E}\widetilde M_i^2(\theta,X_i)\big/I_{n}\to 0,
\nonumber\\ \label{p1-j29}
\sum\limits_{i=1}^n\overline h_i\big({\bf E}|\theta_n^*-{\bf E}_{\neq i}\theta_{n}^*|^2\big)^{p/2}\big({\bf E}\widetilde M_i^2(\theta,X_i)\big)^{1/2}\big/\sqrt{I_{n}}\to 0,
  \end{gather}
где через~ ${\bf E}_{\neq i}$   обозначено условное
математическое ожидание,
  взятое при условии, что   при всех $j\neq i$ фиксированы значения
     $X_{j},$ $j=1,\ldots,n$.
 К вопросу о проверке  условия (\ref{p1-j29})  стоит отметить, что в простейших регулярных случаях  ${\bf E}|\theta_n^*-{\bf E}_{\neq i}\theta_{n}^*|^2$ есть величина порядка $1/n^2$.
   \hfill$\square$
\end{rem}

\section{Обобщения некоторых результатов на случай многомерного параметра }

В этом параграфе часть  результатов \S2 будут обобщены на случай многомерного параметра $\tha\subset \mathbb{R}^m.$

Условимся об обозначениях. Далее  вектора --- это вектор-столбцы размерности $m$, которые будут обозначаться полужирными буквами, а $m$-мерные квадратные матрицы будут обозначаться прямыми  заглавными буквами. Символ $\;^\top$ обозначает транспонирование вектора или матрицы.
 В случае симметричных неотрицательно определенных матриц ${\rm D}$ через  ${\rm D}^{1/2}$ будем обозначать единственную симметричную неотрицательно определенную матрицу, которая удовлетворяет равенству
$\big({\rm D^{1/2}}\big)^2={\rm D}={\rm D}^{1/2}{\rm D}^{1/2\top}. $

 Для
любого вектора ${\bf x}$ символ $\|{\bf x}\|$ обозначает евклидову
норму этого вектора, а через $\|{\rm B}\|^2$ обозначается  сумма квадратов элементов   матрицы ${\rm B}$. Отметим, что таким образом введенная матричная норма обладает всеми свойствами нормы на линейном пространства матриц в силу очевидной изометрии этого пространства с евклидовым пространством размерности $m^2$. Кроме того, для любых
матриц ${\rm B}$ и ${\rm C}$ и любого вектора $\bf x$ справедливы соотношения $\|{\rm BC}\|\leq \|{\rm B}\|\|{\rm C}\|  $ и $\|{\rm B}{\bf x}\|\leq \|{\rm B}\|\|{\bf  x}\|. $
Сходимость векторов или матриц  понимается или
как покоординатную сходимость, или, что эквивалентно,
 как сходимость в соответствующей евклидовой норме.
Символ ${\cal N}({\bf 0},{\rm I})$ обозначает $m$-мерное стандартное нормальное распределение,  а символ $\mathbb{E}$ --- математическое ожидание, взятое по истинному распределению.

В случае $m$-мерного параметра  $\tha\subset \mathbb{R}^m$ определяющие $M$-оценку функции ${\bf M}_i({\bf t},X_i)$ являются $m$-мерными векторами и $M$-оценка $\widetilde\tha_n$ определяется как решение $m$ уравнений, которые в векторной форме имеют следующий вид:
\begin{equation*}   
\sum\limits_{i=1}^n{\bf M}_i({\bf t},X_i)=0.
\end{equation*}

Чтобы в этом случае определить одношаговые $M$-оценки   и сформулировать утверждения об их свойствах, приведем
соответствующие обобщения  условий $(A_1)-(A_5),$ составляющих предположения  $(A).$

 $ ({\bf A}_1).$
Наблюдается выборка объема $n$,  состоящая из независимых элементов   $X_{1},$ $X_{2},\ldots, X_{n}$  со значениями в
 произвольном измеримом пространстве
${\cal X}$ и  распределениями~${\cal L}_{1}$ ${\cal L}_{2},\ldots,{\cal L}_{n}$, зависящими  от    интересующего нас  основного
неизвестного параметра $\tha\in\Tha\subset
\mathbb{R}^m$, где $\Tha$ --- открытое множество.  Кроме того,  эти распределения зависят, вообще говоря, от $n$, и, возможно, зависят  также  от некоторых мешающих  параметров
 произвольной природы.

$ ({\bf A}_2).$  На множестве  $\Tha$ при любом  $i$ заданы зависящие, вообще говоря, от $n$  $m$-мерные вектора   ${\bf M}_i({\bf t},X_i)$ и $m\times m$-мерные матрицы ${\rm M}'_i({\bf t},X_i)$
такие, что
для каждого интервала
$({\bf t}_1,{\bf t}_2)$, целиком лежащего в $\Theta$, с
   вероятностью $1$ имеет место равенство
\begin{equation*}
{\bf M}_i({\bf t}_2,X_i)-
{\bf M}_i({\bf t}_1,X_i)=\int\limits_{0}^{1}
 {\rm M}'_i({\bf t}_1+v({\bf t}_2-{\bf t}_1),X_i)dv\cdot ({\bf t}_2-{\bf t}_1)\qquad \forall ({\bf t}_1,{\bf t}_2)\subset  \Theta,
\end{equation*}
 конечны математические ожидания каждой компоненты матрицы ${\rm M}_i'(\theta,X_i)$
  и справедливы  соотношения
\begin{equation*}  
{\mathbb{ E}}{\bf M}_i(\tha,X_i)={\bf 0}, \qquad {\mathbb E}\|{\bf M}_i(\tha,X_i)\|^2<\infty.
\end{equation*}

\begin{rem}
Отметим, что для почти всех ${\bf t}\in\Tha $
матрица ${\rm M}'_i({\bf t},X_i)$ с вероятностью~$1$ совпадает с матрицей Якоби  вектора ${\bf M}_i({\bf t},X_i)$.
\end{rem}

$ ({\bf A}_3).$ Начиная с некоторого $n$ матрица ${\rm J}_n:=\sum\limits_{i=1}^n{\mathbb E}{\rm M}_i'(\theta,X_i)$ невырождена, а матрица
$
{\rm I}_n:=
\sum\limits_{i=1}^n
{\mathbb E}{\bf M}_i(\tha,X_i){\bf M}_i^\top(\tha,X_i)
$, являющаяся суммой ковариационных матриц
 векторов  $
{\bf M}_i(\tha,X_i),
$ положительно определена.
 \;\;\; Кроме того, \;\;\; $ \|{\rm I}_n^{-1/2}\|\|{\rm J}_n\|\to~\infty,$ \;\;\; $\sup\limits_{n}\|{\rm I}_n^{1/2}\|\|{\rm I}_n^{-1/2}\|<\infty,$
\begin{equation}\label{p1-w5}
\sum\limits_{i=1}^n {\rm M}'_i(\tha,X_i){\rm J}_n^{-1}\stackrel{p}{\to}{\rm I}, \qquad
{\rm I}_n^{-1/2}\sum\limits_{i=1}^n {\bf M}_i(\theta,X_i)\Longrightarrow {\cal N}({\bf 0},{\rm I}).
\end{equation}

Положим $\quad \overline\tau_n(\|\ddelta\|)=\displaystyle\sum\limits_{i=1}^n{\mathbb E}\tau_{i}(\|\ddelta\|,X_i),$ где
 \begin{equation}   \label{p1-wtau}
\tau_{i}(\|\ddelta\|,X_i):=
\begin{cases}
\sup\limits_{{\bf t}
:\; \|{\bf t}-\tha\|\leq\|\ddelta\|}\big\|
\big({\rm M}'_i({\bf t},X_i)-{\rm M}'_i(\tha,X_i)\big){\rm J}_n^{-1} \big\|,\quad \mbox{если }\; [\tha-\ddelta,\tha+\ddelta]\subset\Tha,
\\
\; \infty,\quad \mbox{иначе}.
\end{cases}
\end{equation}

$ ({\bf A}_4).$   Имеет место соотношение
$$\lim\limits_{\|\ddelta\|\to 0}\limsup\limits_{n\to\infty} \overline\tau_n(\|\ddelta\|)= 0.$$

$ ({\bf A}_5).$ Имеется  оценка  $\tha_n^*=\tha_n^*(X_1,\ldots,X_n)$ такая, что
 \begin{gather*}  
  \|{\bf I}_n^{-1/2}\|\|{\bf J}_n\|\|\tha_n^*-\tha\|\overline\tau_n(\|\tha_n^*-\tha\|)\stackrel{p}{\to} 0.
  \end{gather*}

При выполнении условий  $(A) $   одношаговую $M$-оценку  $\tha_{n,M}^{**}$ для параметра $\tha$ определим соотношением
 \begin{gather} \label{p1-w7}
 \tha_{n,M}^{**}=\tha_n^*-
 \left(\sum\limits_{i=1}^n{\rm M}'_i(\tha_n^*,X_i)\right)^{-1}{\sum\limits_{i=1}^n {\bf M}_i(\tha_n^*,X_i)}.
 \end{gather}

 Имеет место следующее обобщение теоремы \ref{p1-t1}.

\begin{theor}                                     \label{p1-wt1}
Пусть  выполнены условия    $(A)$.
 Тогда
оценка $\tha_{n,M}^{**}$  определена с вероятностью, стремящейся к $1$ и
\begin{gather*} 
{\rm I}_n^{-1/2}{\rm J}_n\big(\tha_{n,M}^{**}-\tha\big)\Longrightarrow {\cal N}({\bf 0},{\rm I}).
 \end{gather*}
\end{theor}

В случае многомерного параметра аналог статистики (\ref{p1-777}) определяется соотношением
 \begin{gather} \label{p1-w25}
 \tha_{n,J}^{**}=\tha_n^*-
 {\rm J}_n^{*-1}{\sum\limits_{i=1}^n {\bf M}_i(\tha_n^*,X_i)},.
 \end{gather}
 где ${\rm J}_n^*$ --- некоторая статистика такая, что ${\rm J}_n^*{\rm J}_n^{-1}\stackrel{p}{\to}{\rm I}.$
Имеет место
\begin{theor}                                     \label{p1-wt3}
Пусть  выполнены условия   $(A)$ и
 \begin{gather} \label{p1-w26}
 \|{\rm I}_n^{-1/2}\|\|\tha_n^*-\tha\|\|{\rm J}_n^*-{\rm J}_n\|\stackrel{p}{\to}0,\qquad {\rm J}_n^*{\rm J}_n^{-1}\stackrel{p}{\to}{\rm I}.
 \end{gather}
 Тогда
оценка $\tha_{n,J}^{**}$  определена с вероятностью, стремящейся к $1$ и
\begin{gather*}                         
{\rm I}_n^{-1/2}{\rm J}_n\big(\tha_{n,J}^{**}-\tha\big)\Longrightarrow {\cal N}({\bf 0},{\rm I}).
 \end{gather*}
\end{theor}

\section{Доказательства}

{\bf 5.1.}  Определение  оценки  $\theta_{n,M}^{**}$ из (\ref{p1-7})  перепишем в виде
\begin{equation}                                                                       \label{p1-100}
\theta_{n,M}^{**}-\theta=-\frac{M_n(\theta_n^*)-(\theta_n^*-\theta)M'_n(\theta_n^*)}{M'_n(\theta_n^*)}
\quad\mbox{при}\quad M_n'(t)=\sum\limits_{i=1}^n M_i'(t,X_i).
\end{equation}
Положим $u_{n}^*=\theta_n^*-\theta,$
\begin{equation}  \label{p1-101}
 \rho_{n}(u):=M_n'(\theta+u)-M_n'(\theta),\quad
\overline\rho_{n}(u):=M_n(\theta+u)-M_n(\theta)-uM_n'(\theta).
 \end{equation}
В этом случае равенство из  (\ref{p1-100}) и условия из (\ref{p1-5})
 можно переписать соответственно  в виде
\begin{gather}                                                                     \label{p1-102}
\frac{J_n}{\sqrt{I_n}}(\theta_{n,M}^{**}-\theta)= -\frac{M_n(\theta)/\sqrt{I_n}+\overline\rho_{n}(u_n^*)/\sqrt{I_n}-u_n^*\rho_{n}(u_n^*)/\sqrt{I_n}}
{M_n'(\theta)/J_n+\rho_{n}(u_n^*)/J_n},
  \end{gather}
\begin{gather}                                                                     \label{p1-103}
M_n'(\theta)/J_n\stackrel{p}{\to}1,\qquad M_n(\theta)/\sqrt{I_n}\Longrightarrow {\cal N}(0,1).
  \end{gather}
Представление (\ref{p1-102})  является основным в доказательствах теорем~\ref{p1-t1}, \ref{p1-t2} и \ref{p1-t1w+}.
Нам  также потребуются   соотношения
 \begin{equation}\label{p1-119}
\overline{\rho}_n(u)=
u\int\limits_0^1\rho_n(uv)dv,
\end{equation}
 \begin{equation}\label{p1-119+}
 {\bf E}\gamma_{n}(\delta)=\overline\tau_n(\delta)\quad \mbox{при}\quad
 \gamma_{n}(\delta):=|J_n|^{-1}\sum\limits_{i=1}^n \tau_{i}(\delta,X_i).
\end{equation}
Отметим, что тождество (\ref{p1-119}) немедленно следует из формулы Тейлора с остаточным членом в интегральной форме, поскольку с учетом определений из   (\ref{p1-101})
$$
M_n(\theta+u)=M_n(\theta)+u\int\limits_0^1M_n'(\theta+uv)dv\pm uM_n'(\theta)\equiv M_n(\theta)+uM_n'(\theta)+u\int\limits_0^1\rho_n(uv)dv.
$$

Приведем несколько вспомогательных утверждений, которые  нужны для оценки погрешностей при выводе  теорем \ref{p1-t1}- \ref{p1-t3}.

\begin{lem}\footnote{Эта лемма доказана в [\ref{2222}]. Здесь мы приводим иное более короткое доказательство этого утверждения  по сравнению  с оригинальным.}
 \label{p1-lem5}
Пусть $\alpha_n(\cdot): [0,\infty)\to [0,\infty)$ -- монотонно неубывающий случайный процесс, $\beta_n(\cdot)={\bf E}\alpha_n(\cdot)$ и
 $\beta_n(\eta_n)\stackrel{p}{\to}0$ для некоторой случайной величины  $\eta_n\geq 0$.
 Тогда  $\alpha_n(\eta_n)\stackrel{p}{\to}0$.
\end{lem}
 Д~о~к~а~з~а~т~е~л~ь~с~т~в~о.
 Пусть $\delta$ --- произвольное положительное число. Обозначим $c_\delta=\inf\{t: \beta_n(t)\geq \delta\}$.
 Рассуждая от противного, нетрудно видеть, что $\beta_n(t)\to 0$ при $t\to 0$, а потому $c_\delta\to 0$ при $\delta\to 0.$
 Далее, используя неравенство $
 {\bf P}(\eta_n>c_\delta)\leq {\bf P}(\beta_n(\eta_n)\geq \delta),
 $
 для любого $\varepsilon>0$ получаем оценку
 $$
 {\bf P}(\alpha_n(\eta_n)>\varepsilon)\leq {\bf P}(\alpha_n(c_\delta)>\varepsilon)+{\bf P}(\eta_n>c_\delta)\leq \beta_n(c_\delta)/\varepsilon
+{\bf P}(\beta_n(\eta_n)\geq \delta).
 $$
Это соотношение  доказывает утверждение леммы, поскольку для любого $\varepsilon>0$ выбором $\delta>0$ можно гарантировать сколь угодную малость $\beta_n(c_\delta)/\varepsilon$. \hfill$\square$

\begin{lem} \label{p1-lem1}
Если  выполнены условия $(A)$, то
\begin{equation}  \label{p1-110-}
\rho_{n}(u_n^*)/J_n\stackrel{p}{\to}0,\quad u_n^*\rho_n(u_n^*)/\sqrt{I_n}\stackrel{p}{\to}0\quad\mbox{и}\quad
   \overline\rho_{n}(u_n^*)/\sqrt{I_n}
\stackrel{p}{\to}0,
\end{equation}
а оценка $\theta_{n,M}^{**}$ определена с вероятностью, стремящейся к $1$.
\end{lem}
Д~о~к~а~з~а~т~е~л~ь~с~т~в~о леммы \ref{p1-lem1}. В условиях леммы справедливы сходимости
 \begin{equation}\label{p1-119-}
 \overline \tau_n(|u_n^*|)
 \stackrel{p}{\to }0\qquad\mbox{и}\qquad  |J_n||u_n^*|\overline\tau_n(|u_n^*|)/\sqrt{I_n}\stackrel{p}{\to }0.
 \end{equation}
Действительно, вторая сходимость в (\ref{p1-119-}) совпадает с условием $(A_5)$. Кроме того, с учетом этой сходимости для  любого $\varepsilon>0$
 \begin{eqnarray}\label{p1-118}
 {\bf P}\Big(\overline \tau_n(|u_n^*|)>\varepsilon\Big)= {\bf P}\Big(\overline \tau_n(|u_n^*|)>\varepsilon,\; \frac{|J_n|}{\sqrt{I_n}}|u_n^*|\geq 1\Big)+
 {\bf P}\Big(\overline \tau_n(|u_n^*|)>\varepsilon,\; \frac{|J_n|}{\sqrt{I_n}}|u_n^*|< 1\Big)\leq
 \nonumber\\ \leq
 {\bf P}\Big(|J_n||u_n^*|\overline \tau_n(|u_n^*|)/\sqrt{I_n}>\varepsilon\big)+{\bf P}\Big(\overline \tau_n(|u_n^*|)>\varepsilon,\; |u_n^*|< \frac{\sqrt{I_n}}{|J_n|}\Big)\to 0.\qquad
 \end{eqnarray}
При выводе сходимости в (\ref{p1-118}) мы также учли, что в силу монотонности функции $\overline \tau_n(\cdot)$ и условий $(A_3)$ и $(A_4)$
начиная с некоторых $n$ в правой  части (\ref{p1-118}) под знаком вероятности --- невозможное событие.

Ввиду (\ref{p1-tau}), (\ref{p1-101}), (\ref{p1-119}) и (\ref{p1-119+}), имеем
\begin{equation*} 
\frac{|{\rho}_n(u)|}{|J_n|}\leq \gamma_n(|u|), \quad \frac{|u||\rho_n(u)|}{\sqrt{I_n}} \leq \frac{|J_n||u|\gamma_n(|u|)}{\sqrt{I_n}},\quad
\frac{|{\overline\rho}_n(u)|}{\sqrt{I_n}} \leq \frac{|J_n||u|\gamma_n(|u|)}{\sqrt{I_n}},
\end{equation*}
при этом  функции $\alpha_n(|u|)=\gamma_n(|u|)$  и  $\alpha_n(|u|)=|J_n||u|\gamma_n(|u|)/\sqrt{I_n}$ 
 удовлетворяют условиям леммы \ref{p1-lem5}. А потому сходимости (\ref{p1-110-}) следуют из леммы   \ref{p1-lem5} при $\eta_n=|u_n^*|$, если еще учесть равенство в
 (\ref{p1-119+}).

 Докажем последнее утверждение леммы.
 Поскольку c вероятностью, стремящейся к $1$, величина $\overline\tau_n(|\theta_n^*-\theta|)$ конечна в силу  (\ref{p1-119-}), то ${\bf P}(\theta_n^*\in \Theta)\to 1.$
 Следовательно, согласно условию $(A_2)$, с вероятностью, стремящейся к $1$ определены величины $M_n(\theta_n^*)$ и $M_n'(\theta_n^*).$ А поскольку еще
  $M_n'(\theta_n^*)=M_n'(\theta)+\rho_n(u_n^*), $ то из условия $(A_3)$ и сходимостей в (\ref{p1-103}), (\ref{p1-110-})   заключаем, что ${\bf P}\big(M_n'(\theta_n^*)=0\big)\to 0$.
 Таким образом, с учетом определения (\ref{p1-7}), статистика  $\theta_{n,M}^{**}$ определена с вероятностью, стремящейся к $1$.~\hfill$\square$

 Д~о~к~а~з~а~т~е~л~ь~с~т~в~о теоремы \ref{p1-t1}. Первое утверждение теоремы доказано в лемме \ref{p1-lem1}.
 Из тождества  (\ref{p1-102}) и определения (\ref{p1-1w}), имеем
\begin{gather}                                                                     \label{p1-8w}
\Delta_{n,M}\equiv  -\frac{\big(1-M_n'(\theta)/J_n-\rho_n(u_n^*)/J_n\big)M_n(\theta)/\sqrt{I_n}+\overline\rho_{n}(u_n^*)/\sqrt{I_n}-u_n^*\rho_{n}(u_n^*)/\sqrt{I_n}}
{M_n'(\theta)/J_n+\rho_{n}(u_n^*)/J_n},
  \end{gather}
  Сходимость (\ref{p1-1w}) 
   следует теперь из   тождества  (\ref{p1-8w}), условий (\ref{p1-103}) и  леммы \ref{p1-lem1}.
    \hfill $\square$

\begin{lem} \label{p1-lem2}
Пусть   выполнены условия $(A)$ и  {\rm (\ref{p1-8})}. Тогда
\begin{equation}  \label{p1-110}
\delta_{n,M}:=
\big(M_{n,2}(\theta_{n,M}^{**})-M_{n,2}(\theta)\big)/\sqrt{I_n}\stackrel{p}{\to}0\quad\mbox{при}\quad M_{n,2}^2(t)=\sum\limits_{i=1}^n M_i^2(t,X_i)
\end{equation}
и статистика $d_{n,M}^*$ определена с вероятностью, стремящейся к $1$.
\end{lem}
Д~о~к~а~з~а~т~е~л~ь~с~т~в~о. Положим $u_n^{**}:=\theta_{n,M}^{**}-\theta$.
Имеют место следующие соотношения:
 \begin{equation}  \label{p1-117-}
J_nu_n^{**}/\sqrt{I_n}=O_p(1),\qquad \gamma_n(|u_n^{**}|)
\stackrel{p}{\to}0,\qquad \sum\limits_{i=1}^n
 [M_i'(\theta,X_i)]^2/J_n^2\stackrel{p}{\to}0.
\end{equation}
Первое из этих соотношений  следует из сходимости $\frac{J_n}{\sqrt{I_n}}u_n^{**}\Rightarrow {\cal N}(0,1)$, доказанной в  теореме \ref{p1-t1}. Второе -- из леммы \ref{p1-lem5} при $\eta_n=|u_n^{**}|$ и
 $\alpha_n(\cdot)=\gamma_n(\cdot)$, если только еще учесть равенство в  (\ref{p1-119+}) и тот факт, что $\overline\tau_n(|u_n^{**}|)\stackrel{p}{\to}0$
 в силу первого соотношения в (\ref{p1-117-}), сходимости $|J_n|/\sqrt{I_n}\to \infty $ из $(A_3)$ и собственно условия $(A_4).$ Третье соотношение в (\ref{p1-117-}) совпадает с первым условием в (\ref{p1-8}).

Заметим теперь, что в силу  определений   (\ref{p1-119+}),  (\ref{p1-110})   и  условия $(A_2)$, справедливы соотношения
\begin{eqnarray*}
\sqrt{I_n}|\delta_{n,M}|\leq \Big(\sum\limits_{i=1}^n\big(M_i(\theta_{n,M}^{**},X_i)-M_i(\theta,X_i)\big)^2\Big)^{1/2}=
\Big(\sum\limits_{i=1}^n\Big(\int\limits_{\theta}^{\theta_{n,M}^{**}}M_i'(t,X_i)dt\Big)^2\Big)^{1/2}\leq\nonumber\\
\sqrt{2}\Big(\sum\limits_{i=1}^n\Big(\int\limits_{\theta}^{\theta_{n,M}^{**}}\big(M_i'(t,X_i)-M_i'(\theta,X_i)\big)dt\Big)^2\Big)^{1/2}+
\sqrt{2}\Big(\sum\limits_{i=1}^n\Big(\int\limits_{\theta}^{\theta_{n,M}^{**}}M_i'(\theta,X_i)dt\Big)^2\Big)^{1/2}\leq\nonumber\\
\sqrt{2}|u_n^{**}|\Big(\sum\limits_{i=1}^n\Big(\sup\limits_{t:|t-\theta|\leq |u_n^{**}|}\big|M_i'(t,X_i)-M_i'(\theta,X_i)\big|\Big)^2\Big)^{1/2}+
\sqrt{2}|u_n^{**}|\Big(\sum\limits_{i=1}^n[M_i'(\theta,X_i)]^2\Big)^{1/2}\leq\nonumber\\
\sqrt{2}|J_n||\delta_n^{**}|\gamma_n(|u_n^{**}|)+\sqrt{2}|u_n^{**}|\Big(\sum\limits_{i=1}^n[M_i'(\theta,X_i)]^2\Big)^{1/2}.
\end{eqnarray*}
Следовательно, используя  (\ref{p1-117-}), получаем
 \begin{equation*}  
|\delta_{n,M}|\leq \sqrt{2}\frac{|J_n|}{\sqrt{I_n}}|u_n^{**}|\gamma_n(|u_n^{**}|)+\sqrt{2}
\frac{|J_n|}{\sqrt{I_n}}|u_n^{**}|\Big(\sum\limits_{i=1}^n[M_i'(\theta,X_i)]^2/J_n^2\Big)^{1/2}\stackrel{p}{\to}0
\end{equation*}
и первое  утверждение  леммы доказано.

 По теореме \ref{p1-t1} статистика $\theta_{n,M}^{**}$ является асимптотически нормальной, поэтому  ${\bf P}(\theta_{n,M}^{**}\in \Theta)\to 1.$ А поскольку еще ${\bf P}(\theta_n^*\in \Theta)\to 1$ (см. доказательство леммы \ref{p1-lem1}),
 то  с вероятностью, стремящейся к $1$ определены величины $M_n'(\theta_n^*)$ и $M_{n,2}(\theta_{n,M}^{**})$. Кроме того,
  $M_{n,2}(\theta_{n,M}^{**})/\sqrt{I_n}=\delta_{n,M}+M_{n,2}(\theta)/\sqrt{I_n}, $ то из второго условия в (\ref{p1-8}) и сходимости (\ref{p1-110}) заключаем, что ${\bf P}\big(M_{n,2}(\theta_{n,M}^{**})=0\big)\to 0$.
 Таким образом, с учетом определения (\ref{p1-9}), статистика  $d_{n,M}^{*}$ определена с вероятностью, стремящейся к $1$.
\hfill $\square$

   Д~о~к~а~з~а~т~е~л~ь~с~т~в~о теоремы \ref{p1-t2}. Используя
  тождества в (\ref{p1-100}) и (\ref{p1-102}), а также обозначения
   из (\ref{p1-110}), получаем следующее представление:
\begin{equation}  \label{p1-113}
d_{n,M}^*(\theta_{n,M}^{**}-\theta)=-\frac{M_n(\theta)/\sqrt{I_n}+\overline\rho_{n}
(u_n^*)/\sqrt{I_n} -u_n^*\rho_{n}(u_n^*)/\sqrt{I_n}}
{M_{n,2}(\theta)/\sqrt{I_n}+\delta_{n,M}}.
\end{equation}
Но $M_{n,2}^2(\theta)/I_n^2\stackrel{p}{\to}1 $ в силу второго условия в (\ref{p1-8}).
Утверждения теоремы  вытекают теперь из  лемм \ref{p1-lem1}, \ref{p1-lem2},
сходимостей из (\ref{p1-103}) и тождества (\ref{p1-113}).
 \hfill $\square$

\begin{lem} \label{p2-lem3}
Пусть выполнены условия теоремы {\rm\ref{p1-t1w+}}.   Тогда
\begin{equation}  \label{p2-111}
\delta_{n,\rho}:=
\Big|\frac{u_n^*\rho_n(u_n^*)}{\sqrt{I_n}}-
\frac{\overline\rho_{n}(u_n^*)}{\sqrt{I_n}}-\frac{J_n}{\sqrt{I_n}}(u_n^*)^2\frac{\sum\limits_{i=1}^n{\bf
E}M_i''(\theta,X_i)}{2J_n}\Big|=\frac{J_n}{\sqrt{I_n}}(u_n^*)^2 o_p(1).
\end{equation}
\end{lem}
Д~о~к~а~з~а~т~е~л~ь~с~т~в~о леммы   \ref{p2-lem3}. Положим
\begin{equation}\label{p2-121}
\mu_i(u,X_i):=
\begin{cases}\big(M'_i(\theta+u,X_i)-M'_i(\theta,X_i)\big)/u-M''_i(\theta,X_i),\;
u\neq 0,
\\
0, \; u=0,
\end{cases}
\end{equation}
\begin{equation}\label{p2-125}
\overline{\mu}_n(\delta):=\sum_{i=1}^n \sup\limits_{u:|u|\leq\delta, \theta+u\in\Theta} |\mu_i(u,X_i)|/|J_n|.
\end{equation}
Прежде всего покажем, что
\begin{equation}                                                                 \label{p1-128}
\overline{\mu}_n(u_n^*)\stackrel{p}{\to} 0.
\end{equation}
Действительно,  поскольку
$
\mu_i(u,X_i)=\displaystyle\int\limits_0^1 M_i''(\theta+uv,
X_i)dv-M_i''(\theta,X_i),
$
то
$$
\overline\mu_n(\delta)
\leq \sum\limits_{i=1}^n
 \sup\limits_{t\in\Theta:|t-\theta|\leq\delta} |M_i''(t,X_i)-M_i''(\theta,X_i)|/|J_n|.
 $$
 Следовательно, при всех
неслучайных $\varepsilon>0$, $\alpha>0$ и некотором $\delta=\delta(\alpha)>0$, имеем
\begin{eqnarray}\label{p2-130}
\mathbf{P}\big(\overline{\mu}_n(u_n^*)>\varepsilon\big)\leq\mathbf{P}
\big(|u_n^*|>\delta\big)+\mathbf{P}\big(\overline{\mu}_n(\delta)>\varepsilon\big), \nonumber\\
\limsup\limits_{n\to \infty}\mathbf{P}\big(\overline{\mu}_n(u_n^*)>\varepsilon\big)
\leq\limsup\limits_{n\to \infty}\mathbf{P}\big(\overline{\mu}_n(\delta)>\varepsilon\big)= P(\delta)<\alpha.
\end{eqnarray}
При выводе (\ref{p2-130}), доказывающем (\ref{p1-128}), мы учли, что
 $P(\delta)\to 0$ при $\delta\to 0$ ввиду~$(A_6)$.

Завершим доказательство леммы. Сравнивая определения (\ref{p1-101}),
(\ref{p2-121}) и определение $\overline{M''_n}(\theta)$ из $(A_6)$, имеем
$$
\rho_n(vu_n^*)/J_n-vu_n^*\sum\limits_{i=1}^n{\bf E}M_i''(\theta,X_i)/J_n=vu_n^*\overline{M''_n}(\theta)+vu_n^*\sum_{i=1}^n\mu_i(vu_n^*,X_i)/J_n
.
$$
Из этого равенства и (\ref{p2-125}), c учетом условия $(A_6)$ и сходимости (\ref{p1-128}),  получаем
\begin{eqnarray}                                                                     \label{p2-128}
\sup\limits_{v\in[0,1]}\big|\rho_n(vu_n^*)/J_n-vu_n^*\sum\limits_{i=1}^n{\bf E}M_i''(\theta,X_i)/J_n\big|\leq \nonumber\\ |u_n^*|\big(|\overline{M_n''}(\theta)|+\overline{\mu}_n(u_n^*)\big)=|u_n^*|o_p(1).
\end{eqnarray}
 Используя теперь  определение (\ref{p2-111}) и представление (\ref{p1-119}), имеем
\begin{eqnarray*}
\delta_{n,\rho}=\Big|\frac{J_n}{\sqrt{I_n}}u_n^*\Big(\frac{\rho_n(u_n^*)}{J_n}-u_n^*\frac{\sum\limits_{i=1}^n{\bf E}M_i''(\theta,X_1)}{J_n}\Big) -\nonumber\\
\frac{J_n}{\sqrt{I_n}}u_n^*\int_0^1\Big(\frac{\rho_n(vu_n^*)}{J_n}-vu_n^*\frac{\sum\limits_{i=1}^n{\bf E}M_i''(\theta,X_1)}{J_n}\Big)dv
\Big|.
\end{eqnarray*}
Следовательно, с учетом (\ref{p2-128}),
$$
\delta_{n,\rho}\leq 2\frac{|J_n|}{\sqrt{I_n}}|u_n^*|\sup_{v\in[0,1]}
\Big|\rho_n(vu_n^*)/J_n-vu_n^*\sum\limits_{i=1}^n{\bf E}M_i''(\theta,X_i)/J_n\Big|=\frac{J_n}{\sqrt{I_n}}(u_n^*)^2 o_p(1)
$$
и утверждение леммы  \ref{p2-lem3} установлено.
  \hfill $\square$

 Д~о~к~а~з~а~т~е~л~ь~с~т~в~о теоремы \ref{p1-t1w+}. Из тождества  (\ref{p1-8w}),  леммы \ref{p1-lem1} и условий (\ref{p1-103}),  имеем
\begin{equation}  \label{p2-114w}
\Delta_{n,M}=\frac{o_p(1)+u_n^*\rho_{n}(u_n^*)/\sqrt{I_n}-
\overline\rho_{n}(u_n^*)/\sqrt{I_n}} {1+o_p(1)}.
\end{equation}
Соотношение  (\ref{p1-3w}) следует теперь из
(\ref{p2-114w}) и леммы \ref{p2-lem3}. Второе утверждение теоремы очевидно.
 \hfill $\square$

\begin{lem} \label{p1-lem2h-}
Пусть   выполнены условия $(A)$,  {\rm (\ref{p1-8})}~и~{\rm (\ref{p1-555-})}. Тогда
 статистика $d_{n,J}^*$ определена с вероятностью, стремящейся к $1$ и
\begin{equation}  \label{p1-110h-}
\delta_{n,J}:=
\big(M_{n,2}(\theta_{n,J}^{**})-M_{n,2}(\theta)\big)/\sqrt{I_n}\stackrel{p}{\to}0.
\end{equation}
\end{lem}
Д~о~к~а~з~а~т~е~л~ь~с~т~в~о этого утверждения с очевидными изменениями повторяет вывод леммы \ref{p1-lem2}.

   Д~о~к~а~з~а~т~е~л~ь~с~т~в~о теоремы \ref{p1-t1!-!!}. Из определений (\ref{p1-777}), (\ref{p1-101}) и обозначения, введенного в теореме,  получаем следующее тождество:
   \begin{gather}                                                                     \label{p1-102-}
\frac{J_n}{\sqrt{I_n}}(\theta_{n,J}^{**}-\theta)= \frac{-M_n(\theta)/\sqrt{I_n}-\overline\rho_{n}(u_n^*)/\sqrt{I_n}+u_n^*(J_n^*-J_n)/\sqrt{I_n}-J_nu_n^*\overline{M_n'}(\theta)/\sqrt{I_n}}
{J_n^*/J_n}.
  \end{gather}
Сходимость (\ref{p1-24g-}) следует теперь из тождества (\ref{p1-102-}), леммы \ref{p1-lem1} и условия (\ref{p1-555-}). Доказательство
 того, что $\theta_{n,J}^{**}$ определена с вероятностью, стремящейся к $1$, с очевидными изменениями повторяет аналогичные рассуждения  об оценке
 $\theta_{n,M}^{**}$ из леммы \ref{p1-lem1}.

 Далее, из равенства (\ref{p1-102-}) и определений (\ref{p1-9g-}), (\ref{p1-110}) и (\ref{p1-110h-}), имеем
  \begin{gather}                                                                     \label{p1-113-}
d_{n,J}^*(\theta_{n,J}^{**}-\theta)= \frac{-M_n(\theta)/\sqrt{I_n}-\overline\rho_{n}(u_n^*)/\sqrt{I_n}+u_n^*(J_n^*-J_n)/\sqrt{I_n}-J_nu_n^*\overline{M_n'}(\theta)/\sqrt{I_n}}
{M_{n,2}(\theta)/\sqrt{I_n}+\delta_{n,J}}.
  \end{gather}
  Доказательство второй части теоремы с очевидными изменениями повторяет вывод теоремы \ref{p1-t2}, если вместо
  равенства (\ref{p1-113}) воспользоваться соотношением (\ref{p1-113-}) и учесть утверждение леммы \ref{p1-lem2h-}.
 \hfill $\square$

  Д~о~к~а~з~а~т~е~л~ь~с~т~в~о утверждений из замечания  \ref{p1-z15}.
  Пусть выполнены ограничения из (\ref{p1-415}). В этом случае для всех достаточно больших $n$
  \begin{eqnarray}\label{p1-267-}
  \frac{|J_n|}{\sqrt{I_n}}|u^*_n|\overline\tau_n(|u_n^*|)\leq \frac{|J_n|}{\sqrt{I_n}}|u^*_n|^{1+\alpha}o_p(1)
=\nonumber\\ \Big(\big(|J_n|/\sqrt{I_n}\big)^{1/(1+\alpha)}|u_n^*|\Big)^{1+\alpha}o_p(1)=
\big(O_p(1)\big)^{1+\alpha}o_p(1) \stackrel{p}{\to}0.
\end{eqnarray}
Если же выполнены  условия из (\ref{p1-416}), то нужно в этих рассуждениях поменять местами символы $o_p(1)$ и $O_p(1)$.
\hfill $\square$

  Д~о~к~а~з~а~т~е~л~ь~с~т~в~о утверждений из замечания  \ref{p1-z15+}.
  Пусть выполнены ограничения из (\ref{p1-415}) и (\ref{p1-415!}). В этом случае
  имеет место (\ref{p1-267-}), т.е.  справедливо $(A_5)$. Кроме того, т.к. $|J_n|/\sqrt{I_n}\to \infty$, то в силу второго условия из (\ref{p1-415}) выполнено $u_n^*\stackrel{p}{\to} 0.$  Следовательно, для всех достаточно больших $n$
 \begin{multline*}
|J_n(\theta^*)/J_n-1|=|u^*_n|^\alpha o_p(1)\stackrel{p}{\to}~0,\\
\frac{|J_n|}{\sqrt{I_n}}|u^*_n|\big|J_n(\theta^*)/J_n-1\big|=\frac{|J_n|}{\sqrt{I_n}}|u^*_n|^{1+\alpha}o_p(1)=\\
=\Big(\big({|J_n|}{\sqrt{I_n}}\big)^{1/(1+\alpha)}|u_n^*|\Big)^{1+\alpha}o_p(1)=
\big(O_p(1)\big)^{1+\alpha}o_p(1)\underset{p}\to 0,
\end{multline*}
\begin{equation*}
\Big(\frac{|J_n|}{\sqrt{I_n}}\Big)|u^*_n|
\overline{M_n'}(\theta)=\Big(\frac{|J_n|}{\sqrt{I_n}}\Big)^{\alpha/(1+\alpha)}|u_n^*|\cdot
\overline{M_n'}(\theta)\Big(\frac{|J_n|}{\sqrt{I_n}}\Big)^{1/(1+\alpha)}=O_p(1)\cdot
o_p(1)\stackrel{p}{\to} 0,
\end{equation*}
т.е. выполнены все условия в (\ref{p1-555-}).

  Чтобы
извлечь (\ref{p1-555-})  в случае, когда выполнены  предположения  (\ref{p1-416}) и (\ref{p1-416!})     --- нужно
повторить приведенные в этом доказательстве рассуждения,  поменяв местами символы $o_p(1)$ и $O_p(1)$.
\hfill $\square$

Положим
$$
\widetilde{\rho}_{n}^M(u):=\sum\limits_{i=1}^nh_i(\theta)\big(\widetilde M_i(\theta+u,X_i)-\widetilde M_i(\theta,X_i)-u \widetilde M_i'(\theta,X_i)\big),
$$
\begin{eqnarray*}
\widetilde{\rho}_{n}^{Mh}(u):=\sum\limits_{i=1}^n\big(h_i(\theta+u)-h_i(\theta)\big)\big(\widetilde M_i(\theta+u,X_i)-\widetilde M_i(\theta,X_i)-u \widetilde  M_i'(\theta,X_i)\big),\\
\widetilde{\rho}_{n}^h(u):=u\sum\limits_{i=1}^n\big(h_i(\theta+u)-h_i(\theta)\big) \widetilde M_i'(\theta,X_i).
\end{eqnarray*}

\begin{lem}\label{p1-lem11w}
Если выполнено условие $(A^-)$, то
\begin{equation*} 
\widetilde{\rho}_{n}^M(u_n^*)/\sqrt{I_{n}}\stackrel{p}{\to}0,\quad
\widetilde{\rho}_{n}^{Mh}(u_n^*)/\sqrt{I_{n}}\stackrel{p}{\to}0\quad\mbox{и}\quad
\widetilde{\rho}_{n}^h(u_n^*)/\sqrt{I_{n}}\stackrel{p}{\to}0.
\end{equation*}
  \end{lem}
Д~о~к~а~з~а~т~е~л~ь~с~т~в~о. Аналогично (\ref{p1-119}), справедливо равенство
 \begin{equation*}
 \widetilde M_i(\theta+u,X_i)- \widetilde M_i(\theta,X_i)-u \widetilde M_i'(\theta,X_i)=u\int\limits_0^1\big( \widetilde M_i'(\theta+uv,X_i)- \widetilde M_i'(\theta,X_i)\big)dv.
\end{equation*}
Следовательно, с учетом условия $(A^-)$,  имеем
\begin{equation}  \label{p1-j50}
\begin{split}
\frac{|\widetilde\rho_n^M(u)|}{\sqrt{I_n}}\leq \frac{|u|^{1+q}\sum\limits_{i=1}^n|h_i(\theta)|\overline{M_i'}}{\sqrt{I_n}}, \qquad
\frac{|\widetilde\rho_n^{Mh}(u)|}{\sqrt{I_n}}\leq \frac{|u|^{1+p+q}\sum\limits_{i=1}^n\overline h_i\overline{M_i'}}{\sqrt{I_n}}, \\
|\widetilde\rho_n^h(u)|/\sqrt{I_n}\leq |u|^{1+p}\sum\limits_{i=1}^n\overline h_i|\widetilde M_i'(\theta,X_i)|/\sqrt{I_n},
\end{split}
\end{equation}
при этом  функции в правой части соотношений (\ref{p1-j50}) удовлетворяют условиям леммы \ref{p1-lem5}, если их рассматривать в качестве $\alpha_n(|u|)$. Утверждения леммы  следуют теперь  из леммы   \ref{p1-lem5} при $\eta_n=|u_n^*|$ и условий $(A^-)$.
\hfill$\square$

Д~о~к~а~з~а~т~е~л~ь~с~т~в~о теоремы \ref{p1-t1!-!!w}. Обозначим через $\delta_n^J$ и  $\delta_n^M$ соответственно  левые части в двух последних условиях в (\ref{p1-17w}), а через $\delta_n^h$ --- левую часть в условии (\ref{p1-19w}).
Тогда из определения (\ref{p1-18w}) оценки $\theta_{n,J}^{**}$ и определений, введенных перед леммой \ref{p1-lem11w}, имеем
\begin{gather*}
\frac{J_n}{\sqrt{I_n}}(\theta_{n,J}^{**}-\theta)= \frac{-M_n(\theta)/\sqrt{I_n} +
\delta_n^J-\delta_n^M-\delta_n^h-\big(\widetilde\rho_{n}^{Mh}(u_n^*)+\widetilde\rho_{n}^{h}(u_n^*)+
\widetilde\rho_{n}^{M}(u_n^*)\big)/\sqrt{I_n}}
{J_n^*/J_n},
  \end{gather*}
  где $M_n(\theta)/\sqrt{I_n}\equiv \sum\limits_{i=1}^nh_i(\theta)\widetilde M_i(\theta,X_i)/\sqrt{I_n}\Rightarrow {\cal N}(0,1).$
Первое утверждение теоремы  следует теперь из этого тождества, условий из  $(A^-)$ и  леммы \ref{p1-lem11w}. Второе утверждение теоремы в силу $(A^-)$  выполнено очевидным образом.
\hfill$\square$

{\bf 5.2.}  Д~о~к~а~з~а~т~е~л~ь~с~т~в~о теоремы \ref{p1-t3}.
Пусть $J_n<0$ и $\delta=\delta_\theta$ таково, что выполнено
  (\ref{p1-11}) и (\ref{p1-11-}). Имеем
 \begin{eqnarray*} 
\frac{\displaystyle\sup\limits_{t:|t-\theta|\leq
\delta_\theta}\sum M'_i(t,X_1)}{|J_n|}\leq\frac{\displaystyle \sum\sup\limits_{t:|t-\theta|\leq
\delta_\theta} M'_i(\theta,X_1)}{|J_n|}=\nonumber\\ =
\frac{\displaystyle \sum\sup\limits_{t:|t-\theta|\leq
\delta_\theta} M'_i(\theta,X_1)}{\displaystyle \sum{\bf E}\sup\limits_{t:|t-\theta|\leq
\delta_\theta} M'_i(\theta,X_1)}\frac{\displaystyle \sum{\bf E}\sup\limits_{t:|t-\theta|\leq
\delta_\theta} M'_i(\theta,X_1)}{|J_n|}\equiv \zeta_n\cdot\frac{\displaystyle \sum{\bf E}\sup\limits_{t:|t-\theta|\leq
\delta_\theta} M'_i(\theta,X_1)}{|J_n|},
\end{eqnarray*}
где через $\zeta_n$
 обозначена левая часть в условии  (\ref{p1-11}), ввиду которого $\zeta_n\stackrel{p}{\to}1$.
При этом, поскольку $\limsup\limits_{n\to\infty}\overline\tau_n(\delta_\theta)<1, $ то существует $\varepsilon>0$ такое, что при достаточно больших $n$  
имеют место оценки
$\overline\tau_n(\delta_\theta)<1-\varepsilon$ и
 \begin{equation*} 
\frac{\displaystyle \sum{\bf E}\sup\limits_{t:|t-\theta|\leq
\delta_\theta}M'_i(t,X_1)}{|J_n|}\leq\frac{\displaystyle \sum{\bf E}M'_i(\theta,X_1)}{|J_n|}+ \overline\tau_n(\delta_\theta)\leq -1+(1-\varepsilon)=-\varepsilon<0.
\end{equation*}
Наконец, с учетом второй сходимости в (\ref{p1-103}) и условия $\sqrt{I_n}/|J_n|\to 0,$ имеем
$$
\frac{M_n(t)}{|J_n|}=\frac{M_n(\theta)}{\sqrt{I_n}}\frac{\sqrt{I_n}}{|J_n|}+\frac{\displaystyle\sum\int\limits_{\theta}^t
M'_n(z)dz}{|J_n|}=o_p(1)+\frac{\displaystyle \sum \int\limits_{\theta}^t M'_n(z)dz}{|J_n|},
$$
а потому (с учетом что функция $M_n'(t)$
 отделена от нуля)  $M_n(\theta-\delta_\theta)>0$ и $M_n(\theta+\delta_\theta)~<~0.$
Следовательно,
c вероятностью, стремящейся к 1,
строго монотонная (убывающая) на интервале $(\theta-\delta_\theta,\theta+\delta_\theta)$  функция $M_n(t)$
меняет знак на концах интервала, тем самым на этом интервале существует
 единственный нуль  $\widetilde\theta_n(\theta)$ этой   функции.
Случай $J_n>0$  рассматривается аналогично.

Из   равенства
$0=M_n\big(\widetilde\theta_n(\theta)\big)=M_n(\theta)+\displaystyle\int\nolimits_{\theta}^{\widetilde\theta_n(\theta)}M_n'(t)dt$ и (\ref{p1-119}), имеем
\begin{equation}\label{p1-148}
\frac{J_n}{\sqrt{I_n}}\big(\widetilde\theta_n(\theta)-\theta\big)=-\frac{M_n(\theta)}{\sqrt{I_n}}\left({M_n'(\theta)}\big/{J_n}+
{\displaystyle\int\limits_0^1\rho_n\big(v(\widetilde\theta_n(\theta)-\theta) \big)}dv\big/{J_n} \right)^{-1}.
\end{equation}
Так как
$\limsup\limits_{n\to\infty}\overline\tau_n(\delta)\to 0$ при $\delta\to 0$ и условие (\ref{p1-11}) выполнено при всех $\delta$ таких, для которых справедливо соотношение (\ref{p1-11-}),   то
в силу доказанного при  сколь угодно малых $\delta$  на интервале $(\theta-\delta,\theta+\delta)$ c вероятностью, стремящейся к $1$  существует  единственный локальный нуль $\widetilde\theta_n(\theta), $ тем самым
$\widetilde\theta_n(\theta)\stackrel{p}{\to}\theta $ и, как нетрудно видеть,
$\overline\tau_n\big(|\widetilde\theta_n(\theta)-\theta|\big)\stackrel{p}{\to} 0$.
А поскольку еще  $\Big|\displaystyle\int\limits_0^1\rho_n\big(v(\widetilde\theta_n(\theta)-\theta) \big)dv\Big|\leq |J_n|^{-1} \gamma_n\big(|\widetilde\theta_n(\theta)-\theta|\big)$ и ${\bf E}\gamma_n(\delta)=\overline\tau_n(\delta),$ то
 $\displaystyle\int\limits_0^1\rho_n\big(v(\widetilde\theta_n(\theta)-\theta) \big)dv/|J_n|\stackrel{p}{\to}0$ в силу леммы \ref{p1-lem5}.
Эта сходимость вместе с представлением (\ref{p1-148}) и соотношениями (\ref{p1-103}) доказывают~(\ref{p1-38}).~\hfill$\square$

 Д~о~к~а~з~а~т~е~л~ь~с~т~в~о теоремы \ref{p1-co1}.
В силу   известной оценки погрешности приближения интеграла Римана интегральными суммами, имеем
\begin{eqnarray*}
\sum\limits_{i=1}^n\Delta z_{ni}f(\theta, z_{n:i})-\int\limits_c^df(\theta,z)dz= \delta_n-\int\limits_{z_{n:n}}^d f(\theta,z)dz,
\end{eqnarray*}
где $|\delta_n|\leq \sum\limits_{i=1}^n\omega_{f,\theta}(\Delta z_{ni})\Delta z_{ni}.$ Следовательно, с учетом (\ref{p1-m9}), выполнено
\begin{equation}\label{p1-m20}
\alpha_n^{1/p}\Big (\sum\limits_{i=1}^n\Delta z_{ni}f(\theta,z_{n:i}) - T(\theta)\Big)\equiv\alpha_n^{1/p} \Big(\sum\limits_{i=1}^n\Delta z_{ni}f(\theta, z_{n:i})-\int\limits_c^df(\theta,z)dz\Big)\to 0.
\end{equation}
В силу же сходимости (\ref{p1-m9-}) и неравенства Чебышева со вторым моментом, справедливо соотношение
 $\alpha_n^{1/p}\sum\limits_{i=1}^n\Delta z_{ni}\varepsilon_{ni}\stackrel{p}{\to}0$. Для доказательства  $\alpha_n$-состоятельности оценки $\theta_n^*$ остается  воспользоваться следующей цепочкой соотношений:
\begin{eqnarray*}
|\theta_n^*-\theta|\equiv \Big|T^{-1}\Big(\sum\limits_{i=1}^n\Delta z_{ni}X_{ni}\Big)-T^{-1}T(\theta)\Big|\leq \\ K\Big|\sum\limits_{i=1}^n\Delta z_{ni}X_{ni}-T(\theta)\Big|^p\leq K\Big|\sum\limits_{i=1}^n\Delta z_{ni}f(\theta,z_{n:i})-T(\theta)\Big|^p+
K\Big|\sum\limits_{i=1}^n\Delta z_{ni}\varepsilon_{ni}\Big|^p.
\end{eqnarray*}

Докажем второе утверждение теоремы.
В силу сходимости $d_n\to 0$ и неравенства Чебышева со вторым моментом, справедливо соотношение
 $\sum\limits_{i=1}^n\Delta z_{ni}\varepsilon_{ni}\stackrel{p}{\to}0. $ Следовательно,
  с учетом представления  (\ref{p1-m20}) и условий второй части теоремы
\begin{equation}\label{p1-m5+}
\sum\limits_{i=1}^n\Delta z_{ni}X_{ni}=\sum\limits_{i=1}^n\Delta z_{ni}f(\theta,z_{n:i}) +\sum\limits_{i=1}^n\Delta z_{ni}\varepsilon_{ni}\stackrel{p}{\to}T(\theta).
\end{equation}
Далее, в силу формулы конечных приращений справедливо тождество
\begin{eqnarray}\label{p1-m6}
\theta_n^*-\theta=T^{-1}\Big(\sum\limits_{i=1}^n\Delta z_{ni}X_{ni} \Big)-T^{-1}\big(T(\theta)\big)=\nonumber\\
\big[T^{-1}(\widetilde T)\big]'\Big(\sum\limits_{i=1}^n\Delta z_{ni}f(\theta,z_{n:i})-T(\theta)+\sum\limits_{i=1}^n\Delta z_{ni}\varepsilon_{ni}\Big),
\end{eqnarray}
где $\widetilde T$ лежит между $\sum\limits_{i=1}^n\Delta z_{ni}X_{ni} $ и $T(\theta),$ а потому $\widetilde T\stackrel{p}{\to} T(\theta)$, если еще учесть (\ref{p1-m5+}).
Следовательно, $\big[T^{-1}(\widetilde T)\big]'=\big[T'\big(T^{-1}(\widetilde T)\big)\big]^{-1}\stackrel{p}{\to} \big[T'(\theta)\big]^{-1}$. Второе утверждение теоремы нетрудно извлечь теперь из   представления (\ref{p1-m6}).
 \hfill$\square$

{\bf 5.3. } В этом разделе докажем утверждения \S~4.

Положим ${\bf u}_{n}^*=\tha_n^*-\tha,$
\begin{equation}  \label{p1-w20}
\begin{split}
{ \rm R}_{n}({\bf u}):=\sum\limits_{i=1}^n \big({\rm M}_i'(\tha+{\bf u},X_i)-{\rm M}_i'(\tha,X_i)\big){\rm J}_n^{-1},\\
\rrho_{n}({\bf u}):={\rm I}_n^{-1/2}\sum\limits_{i=1}^n\big({\rm M}'_i(\tha+{\bf u},X_i)-{\bf M}'_i(\tha,X_i)\big){\bf u},\\
 \overline\rrho_{n}({\bf u}):={\rm I}_n^{-1/2}\sum\limits_{i=1}^n\big({\bf M}_i(\tha+{\bf u},X_i)-{\bf M}_i(\tha,X_i)-{\rm M}_i'(\tha+{\bf u},X_i){\bf u}\big).
 \end{split}
 \end{equation}

Аналогом леммы \ref{p1-lem1} является следующее утверждение.

\begin{lem}\label{p1-wlem-1}
Пусть выполнены условия $(A).$ Тогда
$$
{ \rm R}_{n}({\bf u}_n^*)\stackrel{p}{\to}{\rm 0}, \qquad \rrho_{n}({\bf u}_n^*)\stackrel{p}{\to}{\bf 0},
\qquad \overline\rrho_{n}({\bf u}_n^*)\stackrel{p}{\to}{\bf 0}
$$
и статистика $\tha_{n,M}^{**}$ определена с вероятностью, стремящейся к $1$.
\end{lem}
Д~о~к~а~з~а~т~е~л~ь~с~т~в~о. Заметим прежде всего, что в силу $(A_2)$
 \begin{equation}\label{p1-w33}
 {\bf M}_i(\tha+{\bf u},X_i)-{\bf M}_i(\tha,X_i)-{\rm M}_i'(\tha,X_i){\bf u}=\int\limits_0^1\big({M}_i'(\tha+v{\bf u})-{\rm M}_i'(\tha,X_i)\big)dv{\bf u}.
 \end{equation}
 Кроме того, при выполнении условий (A) справедливы  сходимости
 \begin{equation}\label{p1-w119-}
 \overline \tau_n(\|{\bf u}_n^*\|)
 \stackrel{p}{\to }0\qquad\mbox{и}\qquad  \|{\rm I}_n^{-1/2}\|\|{\rm J}_n\|\|{\bf u}_n^*\|\overline\tau_n(\|{\bf u}_n^*\|)\stackrel{p}{\to }0.
 \end{equation}
Действительно, вторая сходимость в (\ref{p1-w119-}) совпадает с условием $(A_5)$. Кроме того, с учетом этой сходимости для  любого $\varepsilon>0$
 \begin{eqnarray}\label{p1-w118}
 {\bf P}\Big(\overline \tau_n(\|{\bf u}_n^*\|)>\varepsilon\Big)= {\bf P}\Big(\overline \tau_n(\|{\bf u}_n^*\|)>\varepsilon,\; \|{\rm I}_n^{-1/2}\|
 \|{\rm J}_n\|\|{\bf u}_n^*\|\geq 1\Big)+\nonumber\\
 {\bf P}\Big(\overline \tau_n(\|{\bf u}_n^*\|)>\varepsilon,\; \|{\rm I}_n^{-1/2}\|
 \|{\rm J}_n\|\|{\bf u}_n^*\|< 1\Big)\leq
 \nonumber\\ \leq
 {\bf P}\Big(\|{\rm I}_n^{-1/2}\|\|{\rm J}_n\|\|{\bf u}_n^*\|\overline \tau_n(\|{\bf u}_n^*\|)>\varepsilon\big)+{\bf P}\Big(\overline \tau_n(\|{\bf u}_n^*\|)>\varepsilon,\; \|{\bf u}_n^*\|< \frac{1}{\|{\rm I}_n^{-1/2}\|\|{\rm J}_n\|}\Big)\to 0.\qquad
 \end{eqnarray}
При выводе сходимости в (\ref{p1-w118}) использовано, что в силу монотонности функции $\overline \tau_n(\cdot)$ и условий $(A_3)$ и $(A_4)$
начиная с некоторых $n$ в правой  части (\ref{p1-w118}) под знаком вероятности --- невозможное событие.

Далее, в силу  (\ref{p1-wtau}),  (\ref{p1-w20}) и  (\ref{p1-w33}), имеем
\begin{eqnarray}  \label{p1-w31}
\|{\rm R}_n({\bf u})\|\leq\sum\limits_{i=1}^n\tau_i(\|{\bf u}\|,X_i),\quad
\|\rrho_n({\bf u})\|\leq \|{\rm I}_n^{-1/2}\|{\rm J}_n\|\|{\bf u}\|\sum\limits_{i=1}^n\tau_i(\|{\bf u}\|,X_i),\nonumber\\
\|\overline \rrho_n({\bf u})\|\leq \|{\rm I}_n^{-1/2}\|{\rm J}_n\|\|{\bf u}\|\sum\limits_{i=1}^n\tau_i(\|{\bf u}\|,X_i).
\end{eqnarray}
при этом   функции в правых частях соотношений (\ref{p1-w31})
удовлетворяют условиям леммы \ref{p1-lem5}, если их рассматривать в качестве функций
 $\alpha_n(\|{\bf u}\|)$. Таким образом, все утверждения леммы  следуют из леммы   \ref{p1-lem5} при $\eta_n=\|{\bf u}_n^*\|$, если еще учесть равенство $\sum\limits_{i=1}^n\mathbb{E}\tau_i(\|{\bf u}\|,X_i)=\overline \tau_n(\|{\bf u}\|).$

Тот факт, что статистика $\tha_{n,M}^{**}$ определена с вероятностью, стремящейся к~$1$, доказывается аналогично одномерному случаю.~\hfill$\square$

Д~о~к~а~з~а~т~е~л~ь~с~т~в~о теоремы  \ref{p1-wt1}.
Из определения (\ref{p1-w7}), имеем
\begin{gather*} 
\sum\limits_{i=1}{\rm M}_i'(\tha_n^*,X_i) (\tha_{n,M}^{**}-\tha)=
\sum\limits_{i=1}{\rm M}_i'(\tha_n^*,X_i) (\tha_n^*-\tha)-\sum\limits_{i=1}^n{\bf M}_i(\tha_n^*,X_i).
 \end{gather*}
Следовательно, с учетом обозначений (\ref{p1-w20}),
\begin{gather} \label{p1-w22}
{\rm I}_n^{-1/2}{\rm J}_n (\tha_{n,M}^{**}-\tha)=
\Big({\rm I}_n^{-1/2}\big({\rm R}_n({\bf u}_n^*)+{\rm \overline R}_n\big){\rm I}_n^{1/2}+{\rm I} \Big)^{-1}
\Big(\rrho_n({\bf u}_n^*)-\overline\rrho_n({\bf u}_n^*)-\rrho_n \Big),
 \end{gather}
 где в силу условия (\ref{p1-w5})
\begin{gather} \label{p1-w23}
  {\rm \overline R}_n:=\sum\limits_{i=1}{\rm  M}_i(\tha,X_i){\rm J}_n^{-1}-{\rm I}\stackrel{p}{\to }{\rm 0},\qquad
  \rrho_n:={\rm I}_n^{-1/2}\sum\limits_{i=1}{\bf  M}_i(\tha,X_i)\Longrightarrow{\cal N}({\bf 0},{\rm I}).
\end{gather}
Из (\ref{p1-w23}), первого утверждения леммы \ref{p1-wlem-1} и условия $(A_3)$ получаем, что
\begin{gather*} 
  \big \| {\rm I}_n^{-1/2}\big({\rm R}_n({\bf u}_n^*)+{\rm \overline R}_n\big){\rm I}_n^{1/2}\big\|
  \leq \big(\|{\rm R}_n({\bf u}_n^*)\|+\|{\rm \overline R}_n\|\big)\sup\limits_n\|{\rm I}_n^{-1/2}\|\|{\rm I}_n^{1/2}\|\stackrel{p}{\to }0,
\end{gather*}
тем самым ${\rm I}_n^{-1/2}\big({\rm R}_n({\bf u}_n^*)+{\rm \overline R}_n\big){\rm I}_n^{1/2}\stackrel{p}{\to }{\rm 0}$. Сходимость
${\rm I}_n^{-1/2}{\rm J}_n (\tha_{n,M}^{**}-\tha)\Rightarrow{\cal N}({\bf 0},{\rm I})$ следует теперь из представления
(\ref{p1-w22}), двух последних утверждений леммы  \ref{p1-wlem-1} и второй сходимости в (\ref{p1-w23}).
\hfill$\square$

Д~о~к~а~з~а~т~е~л~ь~с~т~в~о теоремы  \ref{p1-wt3}. C учетом  условия (\ref{p1-w26}),  имеем
\begin{gather} \label{p1-w30}
\begin{split}
\rrho_{n1}:={\rm I}_n^{-1/2}({\rm J}_n^*-{\rm J}_n){\bf u}_n^*\stackrel{p}{\to}{\bf 0},\\ \rrho_{n2}:={\rm I}_n^{-1/2}\Big(\sum\limits_{i=1}^n {\rm M}_i'(\tha,X_i)-{\rm J}_n\Big){\bf u}_n^*\stackrel{p}{\to}{\bf 0},\\
{\rm R}_{n1}:={\rm I}_n^{-1/2}\big({\rm J}_n^*{\rm J}_n^{-1}-{\rm I}\big){\rm I}_n^{1/2}\stackrel{p}{\to}{\rm 0}.
\end{split}
 \end{gather}
Первые два соотношения выполнены очевидным образом, а третья сходимость в    (\ref{p1-w30}) следует из  оценки
$$
\|{\rm R}_{n1}\|\leq \|{\rm I}_n^{-1/2}\|\|\big({\rm J}_n^*{\rm J}_n^{-1}-{\rm I}\big)\|\|{\rm I}_n^{1/2}\|\leq
\|\big({\rm J}_n^*{\rm J}_n^{-1}-{\rm I}\big)\|\sup\limits_n\|{\rm I}_n^{-1/2}\|\|{\rm I}_n^{1/2}\|
\stackrel{p}{\to}0.
$$
Используя теперь определение (\ref{p1-w25}) и обозначения (\ref{p1-w20}) и  (\ref{p1-w30}),  получаем тождество
\begin{gather*} 
{\rm I}_n^{-1/2}{\rm J}_n (\tha_{n,J}^{**}-\tha)=
\Big({\rm R}_{n1}+{\rm I} \Big)^{-1}
\Big(\rrho_{n1}-\rrho_{n2}-\overline\rrho_n({\bf u}_n^*)-\rrho_n \Big).
 \end{gather*}
Для установления сходимости ${\rm I}_n^{-1/2}{\rm J}_n (\tha_{n,J}^{**}-\tha)\Rightarrow{\cal N}({\bf 0},{\rm I})$ остается теперь воспользоваться леммой   \ref{p1-wlem-1} и учесть (\ref{p1-w30}) и (\ref{p1-w23}).\hfill$\square$

\vspace{15mm}
\centerline{\Large \textbf{Литература}}
\addcontentsline{toc}{section}{Литература}

\begin{enumerate}

\item\label{2007-B}
{\em Боровков А.А. } Математическая статистика. Москва: Физматлит,
 2007.

\item\label{1981-Dem} {\em Демиденко Е.З.} Линейная и нелинейная регрессия.  М.: Финансы и статистика, 1981.

\item\label{1989-D} {\em Демиденко Е.З.} Оптимизация и регрессия. М.: Наука, 1989.

\item\label{1987-D2} {\em Дрейпер Н., Смит Г.} Прикладной регрессионный анализ. Кн. 2.  М.: Финансы и статистика. 1987.

\item\label{2013-E}  {\em  Ермоленко К.В., Саханенко А.И. } Явные асимптотически нормальные оценки неизвестного  параметра частично-линейной регрессии. ---
 Сиб. электрон. матем. изв., 2013, Т.10, 719-726.

\item\label{1975-Z}
{\em Закс Ш.} Теория статистических выводов. М.:Мир, 1975.

\item\label{2015-K} {\em  Каленчук A.A., Саханенко А.И.} О существовании явных асимптотически нормальных оценок неизвестного параметра логарифмической регрессии. --- Сиб. электрон. матем. изв., 2015, Т.12, 874-883.

\item\label{1991-L}
{\em Леман Э.}  Теория точечного оценивания. М.: Наука, 1991.

\item\label{2015-1?}
{\em  Линке Ю.Ю.} Об уточнении одношаговых оценок Фишера в случае медленно сходящихся предварительных  оценок. --- Теор. вероятн. и ее примен., 2015, Т.60, вып. 1, 80-98.

\item\label{2015-3?}
{\em  Линке Ю.Ю.}
Асимптотические свойства  одношаговых взвешенных $M$-оценок
 с приложениями к задачам регрессии.  ArXiv: 1505.02725v2

\item\label{2011-1}
{\em  Линке Ю.Ю.} Об асимптотике распределения двухшаговых
статистических оценок. ---  Сибирский Математический Журнал, 2011,
Т.~52. \No~4, 841-860.

\item\label{2222}
{\em  Линке Ю.Ю., Саханенко А.И.}
Об условиях асимптотической нормальности одношаговых М-оценок. ---  Сибирский журнал чистой и прикладной математики (в  печати).

\item\label{2013-2}
{\em  Линке Ю.Ю., Саханенко А.И.} Об асимптотике распределения одного
класса двухшаговых статистических оценок многомерного параметра. ---
Математические труды, 2013, Т.~16. \No~1, 89-120.

\item\label{2014-1}
    {\em Линке Ю.Ю., Саханенко А.И.} Об условиях асимптотической нормальности одношаговых оценок Фишера для однопараметрических семейств распределений. ---   Сиб. электрон. матем. изв., 2014, Т.11, 464-475.

\item\label{2000-1}
{\em Линке Ю.Ю. Саханенко А.И.}  Асимптотически нормальное
оценивание параметра в задаче дробно -- линейной регрессии.
 --- Сиб. Матем. Жур., 2000, Т.~41, \No~1, 150-163.

\item\label{2001-2}
{\em   Линке Ю.Ю. Саханенко А.И.}  Явное асимптотически нормальное
оценивание  параметров уравнения Михаэлиса--Ментен. --- Сиб. Матем. Жур., 2001, Т.~42, \No~3, 610-633.

\item\label{2014-S}  {\em  Савинкина Е.М., Саханенко А.И. } Явные оценки неизвестного параметра в одной задаче степенной регрессии. ---
 Сиб. электрон. матем. изв., 2014, Т.11, 725-733.

\item\label{1984-H} {\em Хьюбер П.} Робастность в статистике.  М.:Мир, 1984.

  \item\label{2011-B}   {\em Bergesioa A.,   Yohaia V.} Projection Estimators for Generalized Linear Models. --- J. Amer. Stat. Assoc., 2011,   V. 106(494), 661-671.

\item\label{1975-B} {\em Bickel P.J.}  One-step Huber Estimates in the Linear Model. ---
J. Amer. Stat. Assoc. 1975, V.~70, 428-434.

\item\label{2000-C} {\em Cai Z., Fan J. and Li R.Z.}  Efficient estimation and inferences for varying-coefficient models. --- J.  Amer. Statist. Assoc. 2000, 95, 888-902.

\item\label{2007-C} {\em Cai J., Fan J., Zhou, H. and Zhou, Y.}  Marginal hazard models with varying-coefficients for multivariate failure time data. ---  Ann.  Statist., 2007, 35, 324-354.

\item\label{2009-D} {\em Dette H., Holland-Letz T.}  A geometric characterization of $c$-optimal designs for heteroscedastic regression.  --- Ann. Statist., 2009, V.~37, 6B, 4088-4103.

\item\label{2014-F} {\em Fan, J., Xue, L., and Zou, H.} Strong oracle optimality of folded concave penalized estimation. ---
 Ann.  Statist., 2014,  42, 819-849.

\item\label{1999-F1}{\em Fan J.,  Chen J.} One-step local quasi-likelihood estimation. --- Journal of Royal Statistical Society B, 1999, 61, 927-943.

\item\label{1999-F} {\em Fan J.  Jiang J.}  Variable bandwidth and one-step local M-estimator. --- Science in China, 2000, (Series A), 43, 65-80.

\item\label{1925-F} {\em Fisher R.A. } Theory of statistical estimation. ---   Proc. Camb.
Phil.Soc., 1925, Soc.22,  700-725.

\item\label{1997-H} {\em Heyde C.C.} Quasi-likelihood and its application: a general approach to optimal parameter estimation. 1997. Springer.

\item\label{1971-H} {\em Hoadley B.}   Asymptotic properties of maximum likelihood estimators for the independent not identically distributed case.  --- Ann. Math. Statist., 1971,  T.~42,  N~6, 1977--1991.

\item\label{1973-H} {\em Huber P.J.}   Robust regression: Asymptotics, conjectures and Monte Carlo. --- Ann. Math. Statist. 1973, 799-821.

\item\label{1985-J} {\em Janssen P., Jureckova J., Veraverbeke N.}  Rate of convergence
of one- and two-step M-estimators with applications to maksimum
likelihood and Pitman estimators. ---  Ann. Stat., 1985,  V.~13, N~3,
 1222-1229.

\item\label{1987-J} {\em Jureckova J.,  Portnoy S.} Asymptotics for one-step M-estimators in regression
with application to combining efficiency and high breakdown point. --- Comm. Statist. Theory Methods, 1987, 16, 2187-2200.

\item\label{1990-J} {\em Jureckova J., Sen P.K.}  Effect of the initial estimator on the
asymptotic behavior of the one-step M-estimator. ---  Ann. Inst.
Statist. Math. 1990, V.~42, N~2,  345-357.

\item\label{2006-J} {\em Jureckova J.,  Picek J. } Robust statistical methods with R.
Chapman and Hall, 2006.

\item\label{2012-J1} {\em  Jureckova J.,  Sen P. K.,  Picek J.}
Methodology in Robust and Nonparametric Statistics. Chapman and Hall, Boca Raton, London 2012.

\item\label{2012-J2} {\em Jureckova J.} Tail-behavior of estimators and of their one-step
versions. --- Journal de la Societe Francaise de Statistique, 2012, V. 153  (1), 44-51.

\item\label{2004-L1} {\em Liese F., Vajda I. } A general asymptotic theory of $M$-estimators. I. ---
Mathematical Methods of Statistics. 2004, V.12 (4), 454-477.

\item\label{2004-L2} {\em Liese F., Vajda I. } A general asymptotic theory of $M$-estimators. II. ---
Mathematical Methods of Statistics. 2004, V.13 (1), 82-95.

\item\label{1956-L} {\em Le Cam L. } On the asymptotic theory of estimation and testing hypotheses. ---
 Proceedings of the Third Berkeley Symposium on mathematical statistics and probability, 1956.

\item\label{2011-M} {\em
 Michaelis L., Menten M. L., Johnson K. A., Goody R. S.} The original Michaelis constant: translation of the 1913 Michaelis-Menten paper. ---  Biochemistry.  2011.  V. 50(39),  8264-8269.

\item\label{1989-M} {\em McCullagh P., NelderMuller J.}  Generalized Linear Models. Chapman and Hall. 1989.

\item\label{1994-M1} {\em Muller, Ch.H.}  One-step-M-estimators in conditionally contaminated linear models. ---  Stat. Decis., 1994,  12, 331-342.

\item\label{1994-M2} {\em Muller, Ch.H.}
 Asymptotic behaviour of one-step-M-estimators in contaminated non-linear models. In Asymptotic Statistics, eds. P. Mandl, M. Huskova, Physica-Verlag, Heidelberg, 1994,  395-404.

\item\label{1992-O} {\em Osborne M.R.}  Fisher's method of scoring.  --- International Statistical Review, 1992, V.60 (1), 99-117.

\item\label{1985-R} {\em Reeds J.}  Asymptotic numbers of roots of Cauchy location likelihood equation. --- Ann. Statist., 1985, 13, 775-784.

\item\label{1973-P} {\em Philippou A.N., Roussas G.G.}  Asymptotic distribution of the      likelihood function in the independent not identically distributed case. ---  Ann.  Statist., 1973,  V.~1, N~3, 454-471.

\item\label{1975-P} {\em Philippou A.N., Roussas G.G.}  Asymptotic normality  of the   maximum    likelihood estimate  in the independent not identically distributed case. ---  Ann. Inst.  Statist. Math., 1975,  V.~27, N~1, 45-55.

\item\label{1988-R} {\em Robinson P. M.}  The stochastic difference between econometric statistics. --- Econometrica, 1988, 56, 531-548.

\item\label{2003-S} {\em Seber G.A.F., Wild C.J.} Nonlinear Regression. 2003. John Wiley
and Sons.

\item\label{1980-S} {\em Serfling R.J.} Approximation theorems of mathematical statistics. John Willey and Song,  1980.

\item\label{1992-S} {\em Simpson, D. G., Ruppert, D. and Carroll, R. J.} On One-Step GM Estimates
and Stability of Inferences in Linear Regression. --- Journal of the American
Statistical Association, 1992, 87, 439-450.

\item\label{1986-S} { Stefanski L., Carroll L., Ruppert D.}
Optimally Bounded Score Functions for Generalized Linear Models with Applications to Logistic Regression. ---
Biometrika,  1986, 73 (2): 413-424.

\item\label{2007-V} {\em Verrill S.}  Rate of convergence of k-step Newton estimators to
efficient likelihood estimators. ---  Statistics and Probability
Letters, 2007, 77, 1371--1376.

\item\label{2002-W} {\em  Welsh A.H., Ronchetti E.} A journey in single steps: robust one-step M-estimation in linear regression. ---
Journal of Statistical Planning and Inference, 2002, V. 103  (1-2),  287-310.

\item\label{2004-Y} {\em Yang Y.} Asymptotics of $M$-estimation in non-linear regression.  --- Acta Math. Sinica. 2004, V.~20 (4), 749-760.

\item\label{1979-Y} {\em Yohai V., Maronna R.} Asymptotic behavior of  $M$-estimators for linear model.  --- Ann. Statist., 1979,  V.7 (2), 258-268.

\item\label{2008-Z} {\em Zou H., Li R.} One-step sparse estimates in nonconcave penalized likelihood models.   --- Ann. Statist. 2008, V.~36 (4), 1509-1533.

\end{enumerate}

\end{document}